\documentclass{siamart190516}
\usepackage{amsmath}
\usepackage[T1]{fontenc}
\usepackage[latin9]{inputenc}
\usepackage{array}
\usepackage{verbatim}
\usepackage{float}
\usepackage{units}
\usepackage{mathtools}
\usepackage{bm}
\usepackage{enumitem} 
\usepackage{amssymb}
\usepackage{graphicx}
\usepackage{subcaption}
\usepackage{tcolorbox}
\usepackage{accents}
\usepackage{setspace}
\usepackage{xspace}
\usepackage[hang, flushmargin]{footmisc}

\usepackage[a4paper,
            bindingoffset=0.2in,
            left=2.5cm,
            right=2.5cm,
            top=3.5cm,
            bottom=3.5cm,
            footskip=.3in]{geometry}

\usepackage{algpseudocode}
\usepackage{algorithmicx}

\usepackage{letltxmacro}

\usepackage{thmtools}

\newtheorem{thm}{Theorem}[section]

\newsiamthm{lem}{Lemma}
\newsiamremark{rem}{Remark}

\newlist{thmlist}{enumerate}{1}
\setlist[thmlist]{label=(\roman{thmlisti}), ref=\thethm.(\roman{thmlisti}),noitemsep}
\newlist{lemlist}{enumerate}{1}
\setlist[lemlist]{label=(\roman{lemlisti}), ref=\thelemma.(\roman{lemlisti}),noitemsep}

\usepackage{cleveref}
\usepackage{footnotebackref}
\usepackage{nameref,hyperref}

\Crefname{thm}{Theorem}{Theorems}
\Crefname{thmlisti}{Theorem}{Theorems}

\Crefname{lem}{Lemma}{Lemmas}
\Crefname{lemlisti}{Lemma}{Lemmas}

\crefname{rem}{Remark}{Remarks}
\Crefname{rem}{Remark}{Remarks}

\addtotheorempostheadhook[thm]{\crefalias{thmlisti}{thm}}
\addtotheorempostheadhook[lem]{\crefalias{lemlisti}{lem}}

\crefmultiformat{equation}{Eqs.~(#2#1#3)}%
{ and~(#2#1#3)}{, (#2#1#3)}{ and~(#2#1#3)}
\Crefmultiformat{equation}{Eqs.~(#2#1#3)}%
{ and~(#2#1#3)}{, (#2#1#3)}{ and~(#2#1#3)}

\crefformat{equation}{#2Eq.~(#1)#3}
\Crefformat{equation}{#2Eq.~(#1)#3}

\Crefname{ALC@unique}{Line}{Lines}
\crefname{ALC@unique}{line}{lines}

\usepackage{pgf}
\usepackage{adjustbox}

\begin{document}



\algnewcommand{\IfThenElse}[3]{
	\State \algorithmicif\ #1\ \algorithmicthen\ #2\ \algorithmicelse\ #3}
\algnewcommand\Input{\item[\textbf{Input:}]}%
\algnewcommand\algorithmicinput{\textbf{Input:}}
\algnewcommand\INPUT{\item[\algorithmicinput]}
\algnewcommand\SState{\State \hskip-1.em }
\algnewcommand\algorithmicswitch{\textbf{switch}}
\algnewcommand\algorithmiccase{\textbf{case}}
\algnewcommand\algorithmicassert{\texttt{assert}}
\algnewcommand\Assert[1]{\State \algorithmicassert(#1)}%
\algdef{SE}[SWITCH]{Switch}{EndSwitch}[1]{\algorithmicswitch\ #1\ \algorithmicdo}{\algorithmicend\ \algorithmicswitch}%
\algdef{SE}[CASE]{Case}{EndCase}[1]{\algorithmiccase\ #1}{\algorithmicend\ \algorithmiccase}%
\algtext*{EndSwitch}%
\algtext*{EndCase}%
\def\NoNumber#1{{\def\alglinenumber##1{}\State #1}\addtocounter{ALG@line}{-1}}
\algnewcommand\algorithmicparams{\textbf{Parameters: }}
\algnewcommand\Params{\State \textbf{Parameters:} }%

\definecolor{ccolres}{RGB}{0, 128, 255}
\definecolor{ccolproof}{RGB}{102, 189, 15}
\definecolor{ccoldef}{RGB}{165, 127, 13}

\colorlet{colres}{ccolres!20}
\colorlet{colproof}{ccolproof!20}
\colorlet{coldef}{ccoldef!20}

\newcommand{\noun}[1]{\textsc{#1}}

\newcommand{\sech}{\ensuremath{\operatorname{sech}}}
\newcommand{\GAD}{\ensuremath{\gamma(\alpha + d)}}
\newcommand{\trs}{\noun{trs}}
\newcommand{\btrs}{\noun{btrs}}
\newcommand{\focm}{\noun{focm}}
\newcommand{\lng}{\noun{lng}}
\newcommand{\spdm}{SPD}
\newcommand{\foro}{\noun{foro}}
\newcommand{\rcg}{\noun{rcg}}
\newcommand{\rsd}{\noun{rsd}}
\newcommand{\rgd}{\noun{rgd}\xspace}

\newcommand{\BB}{\mathbb{B}}
\newcommand{\RR}{\mathbb{R}}
\newcommand{\lieig}[2]{\lambda_{#1}\left(#2\right) }
\newcommand{\proj}[1]{\matP_{#1}}
\newcommand{\lmin}{\lambda_{\min}}
\newcommand{\tspM}[2]{T_{#1}{#2}}
\newcommand{\muopt}{\mu_{\star}}
\newcommand{\grad}{\ensuremath{\normalfont{\text{\textbf{grad}}\nobreak\hspace{.05em plus .05em}}}}
\newcommand{\Hess}{\normalfont{\text{\textbf{Hess}}\nobreak\hspace{.05em plus .05em}}}
\newcommand{\lamsub}[1]{\lambda_{#1}}
\newcommand{\rb}{\mathbf{r}}
\newcommand{\ressub}[1]{\rb_{#1}}
\newcommand{\matA}{\ensuremath{{\bm{\mathrm{A}}}}}
\newcommand{\matE}{\ensuremath{{\bm{\mathrm{E}}}}}
\newcommand{\xLN}{\x_{\lng}}
\newcommand{\lLN}{\lambda_{\lng}}
\newcommand{\TNorm}[1]{\|#1\|_{2}}
\newcommand{\TNormS}[1]{\TNorm{#1}^2}
\newcommand{\xk}{\mathbf{x}_{k}}
\newcommand{\xkp}{\mathbf{x}_{k+1}}
\newcommand{\xkt}{\xk^{\T}}
\newcommand{\xkpt}{\xkp^{\T}}
\newcommand{\xki}{\x_{k_{i}}}
\newcommand{\xkip}{\mathbf{x}_{k_{i+1}}}
\newcommand{\xkit}{\xki^{\T}}
\newcommand{\xkipt}{\xkip^{\T}}
\newcommand{\xkj}{\x_{k_{j}}}
\newcommand{\xkjp}{\mathbf{x}_{k_{j+1}}}
\newcommand{\xkjt}{\xkj^{\T}}
\newcommand{\xkjpt}{\xkjp^{\T}}
\newcommand{\ui}{\ensuremath{\u_{i}}}
\newcommand{\uit}{\ensuremath{\u_{i}^{\T}}}
\newcommand{\bx}{\ensuremath{\bar{\x}}}
\newcommand{\by}{\ensuremath{\bar{\y}}}

\newcommand{\vk}{\mathbf{v}_k}
\newcommand{\gk}{\gamma^k}
\newcommand{\lk}{\lamsub{\xk}}
\newcommand{\lkp}{\lamsub{\xkp}}

\newcommand{\calM}{\ensuremath{\mathcal{M}}}

\newcommand{\matM}{\ensuremath{{\bm{\mathrm{M}}}}}
\newcommand{\matMI}{\ensuremath{\matM^{-1}}}
\newcommand{\xMx}{\ensuremath{ (\x^{T} \matMI \x)^{-1} }}
\newcommand{\Mk}{\matM_k}
\newcommand{\MkI}{\Mk^{-1}}
\newcommand{\Mxk}{\Mk \xk}
\newcommand{\MxkI}{ \MkI \xk}
\newcommand{\matH}{\ensuremath{{\bm{\mathrm{H}}}}}
\newcommand{\egrad}[2]{\nabla \bar{#1}(#2)}
\newcommand{\stf}[2]{\ensuremath{{\bm{\mathrm{St}}}(#1, #2)}}
\newcommand{\smm}[1]{\text{Symm}(#1)}
\newcommand{\skk}[1]{\text{Skew}(#1)}

\newcommand{\dist}{\normalfont{\text{dist}\nobreak\hspace{.05em plus .05em}}}
\newcommand{\bzero}{\ensuremath{\mathbf{0}}}

\newcommand{\Asupl}{\matA_{\llng}}
\newcommand{\Alinv}{\matA_{\llng}^{-1}}

\newlength{\dhatheight}
\newcommand{\hhat}[1]{%
	\settoheight{\dhatheight}{\ensuremath{\hat{#1}}}%
	\addtolength{\dhatheight}{-0.35ex}%
	\hat{\vphantom{\rule{1pt}{\dhatheight}}%
		\smash{\hat{#1}}}}

\global\long\def\*#1{\mathbf{#1}}

\global\long\def\xopt{\mathbf{x_{\star}}}

\global\long\def\xlng{\bar{\mathbf{x}}}

\global\long\def\lopt{\mathbf{\lambda_{\star}}}

\global\long\def\llng{\bar{\mathbf{\lambda}}}

\global\long\def\umin{\mathbf{u_{\min}}}

\global\long\def\Umin{\mathbf{U_{\min}}}

\global\long\def\TRS#1#2#3{\text{TRS}_{#1,#2,#3}}

\global\long\def\TRStemplate{\TRS{\*  B}{\*  g}{\Delta}}

\global\long\def\sf#1#2{\mathbf{SF}\left(#1,#2\right)}

\global\long\def\diag{\text{diag}}

\global\long\def\EXP#1#2{\text{EXP}_{#1}\left(#2\right)}

\global\long\def\spec{\text{spec}}

\global\long\def\spn{\text{span}}

\global\long\def\rk#1{\text{rk}\left(#1\right)}

\global\long\def\Okeps{O\left(\nicefrac{k}{\varepsilon}\right)}

\global\long\def\pkeps{\text{\poly\ensuremath{\nicefrac{k}{\varepsilon}}}}

\global\long\def\Ak{\*  A_{k}}

\global\long\def\As{\tilde{\*  A}}

\global\long\def\dag{^{\dagger}}

\global\long\def\psdmat{\mathcal{\boldsymbol{S}}_{n}^{+}}

\global\long\def\psdreal{\mathcal{\boldsymbol{S}}_{1}^{+}}

\global\long\def\symmat{\mathcal{\boldsymbol{S}}_{n}}

\global\long\def\Akr{\*  A_{-k}}

\global\long\def\MM#1{\*{M_{#1} }}

\global\long\def\chik#1{\boldsymbol{\chi}_{k}\left(#1\right)}

\newcommand{\sph}{\ensuremath{\mathbb{S}^{n-1}}}
\newcommand{\sphn}{\mathbb{S}^{(n+1)-1}}
\newcommand{\gsub}[1]{\ensuremath{g_{#1}}}
\newcommand{\tsp}[1]{\ensuremath{T_{#1}\sph}}

\global\long\def\pst{\*{p^{*}}}

\global\long\def\ows{\text{otherwise}}

\global\long\def\lmax{\lambda_{\max}}


\global\long\def\argmin{\text{\ensuremath{\arg\min}}}

\global\long\def\tspopt{\tsp{\pst}}

\global\long\def\projopt{\*{\mathbb{P}_{\pst}}}

\global\long\def\R{\mathbb{R}}

\global\long\def\e{{\mathbf{e}}}

\global\long\def\et#1{{\e(#1)}}

\global\long\def\ef{{\mathbf{\et{\cdot}}}}

\global\long\def\x{{\mathbf{x}}}
\global\long\def\xp{{\mathbf{x}'}}

\global\long\def\w{{\mathbf{w}}}

\global\long\def\xt#1{{\x(#1)}}

\global\long\def\xf{{\mathbf{\xt{\cdot}}}}

\global\long\def\d{{\mathbf{d}}}

\global\long\def\b{{\mathbf{b}}}
\newcommand{\bhat}{\ensuremath{\hat{\b}}}
\newcommand{\bt}{\ensuremath{\b^{\T}}}
\newcommand{\bthat}{\ensuremath{\bhat^{\T}}}

\newcommand{\hatA}{
	\ensuremath{
		\hat{\matA}
	} 
}

\global\long\def\a{{\mathbf{a}}}

\global\long\def\u{{\mathbf{u}}}
\newcommand{\ut}{\ensuremath{{\u}^{\T}}}
\global\long\def\v{{\mathbf{v}}}
\newcommand{\vt}{\ensuremath{{\v}^{\T}}}
\global\long\def\p{{\mathbf{p}}}

\global\long\def\y{{\mathbf{y}}}

\global\long\def\yt#1{{\y(#1)}}

\global\long\def\yf{{\mathbf{\yt{\cdot}}}}

\global\long\def\z{{\mathbf{z}}}
\newcommand{\zt}{\ensuremath{\z^{\T}}}
\global\long\def\v{{\mathbf{v}}}

\global\long\def\e{{\mathbf{e}}}

\global\long\def\h{{\mathbf{h}}}

\global\long\def\s{{\mathbf{s}}}

\global\long\def\g{{\mathbf{g}}}

\global\long\def\c{{\mathbf{c}}}

\global\long\def\p{{\mathbf{p}}}

\global\long\def\f{{\mathbf{f}}}

\global\long\def\xlpr{\x_{\noun{LPR}}}


\global\long\def\rt#1{{\rb(#1)}}

\global\long\def\rf{{\mathbf{\rt{\cdot}}}}

\global\long\def\mat#1{{\ensuremath{\bm{\mathrm{#1}}}}}

\global\long\def\matN{\ensuremath{{\bm{\mathrm{N}}}}}

\global\long\def\matB{\ensuremath{{\bm{\mathrm{B}}}}}

\global\long\def\matC{\ensuremath{{\bm{\mathrm{C}}}}}


\global\long\def\matP{\ensuremath{{\bm{\mathrm{P}}}}}



\newcommand{\matU}{\ensuremath{{\bm{\mathrm{U}}}}}
\newcommand{\matQ}{\ensuremath{{\bm{\mathrm{Q}}}}}
\newcommand{\matT}{\ensuremath{{\bm{\mathrm{T}}}}}
\newcommand{\matV}{\ensuremath{{\bm{\mathrm{V}}}}}
\newcommand{\matD}{\ensuremath{{\bm{\mathrm{D}}}}}
\newcommand{\matI}{\ensuremath{{\bm{\mathrm{I}}}}}

\global\long\def\matR{\mat R}

\global\long\def\matS{\mat S}

\global\long\def\matY{\mat Y}

\global\long\def\matX{\mat X}


\global\long\def\matJ{\mat J}

\global\long\def\matZ{\mat Z}

\global\long\def\matW{\mat W}

\global\long\def\matG{\mat G}

\global\long\def\matL{\mat L}

\global\long\def\S#1{{\mathbb{S}_{N}[#1]}}

\global\long\def\IS#1{{\mathbb{S}_{N}^{-1}[#1]}}

\global\long\def\PN{\mathbb{P}_{N}}

\global\long\def\sgn#1{\text{sign}(#1)}

\global\long\def\InfNorm#1{\|#1\|_{\infty}}

\global\long\def\FNorm#1{\|#1\|_{F}}

\global\long\def\UNorm#1{\|#1\|_{\matU}}

\global\long\def\UNormS#1{\|#1\|_{\matU}^{2}}

\global\long\def\UINormS#1{\|#1\|_{\matU^{-1}}^{2}}

\global\long\def\ANorm#1{\|#1\|_{\matA}}

\global\long\def\BNorm#1{\|#1\|_{\mat B}}

\global\long\def\ANormS#1{\|#1\|_{\matA}^{2}}

\global\long\def\AINormS#1{\|#1\|_{\matA^{-1}}^{2}}

\global\long\def\BINormS#1{\|#1\|_{\matB^{-1}}^{2}}

\global\long\def\BINorm#1{\|#1\|_{\matB^{-1}}}

\global\long\def\ONorm#1#2{\|#1\|_{#2}}

\global\long\def\T{\textsc{T}}

\global\long\def\pinv{\textsc{+}}

\global\long\def\Expect#1{{\mathbb{E}}\left[#1\right]}

\global\long\def\ExpectC#1#2{{\mathbb{E}}_{#1}\left[#2\right]}

\global\long\def\dotprod#1#2{(#1,#2)}

\global\long\def\dotprodM#1#2{(#1,#2)_{\matM}}

\global\long\def\dotprodsqr#1#2{(#1,#2)^{2}}

\global\long\def\Trace#1{{\bf Tr}\left(#1\right)}

\global\long\def\nnz#1{{\bf nnz}\left(#1\right)}

\global\long\def\vol#1{{\bf vol}\left(#1\right)}

\global\long\def\rank#1{{\bf rank}\left(#1\right)}

\global\long\def\diag#1{{\bf diag}\left(#1\right)}

\global\long\def\grad#1{{\bf grad}#1}

\global\long\def\range#1{{\bf range}\left(#1\right)}

\global\long\def\st{\,\,\,\text{s.t.}\,\,\,}

\global\long\def\spacedcolon{\,:\,}

\newcommand{\BS}[1]{{\color{blue}BS: #1}}
\newcommand{\UM}[1]{{\color{olive}UM: #1}}
\newcommand{\HA}[1]{{\color{green}HA: #1}}

\newcommand{\reply}[1]{{\color{blue}#1}}
\newcommand{\DONE}{\reply{Done.}}

	\title{Solving Trust Region Subproblems Using Riemannian Optimization\thanks{This work was supported by Israel Science Foundation grant 1272/17.}}
	\author{Uria Mor\thanks{School of Mathematical Sciences, Tel Aviv University, Tel Aviv, 6997801, Israel (uriamor@mail.tau.ac.il, borisshy@mail.tau.ac.il, haimav@tauex.tau.ac.il).} \and Boris Shustin\footnotemark[2] \and Haim Avron\footnotemark[2]}
	\maketitle

	\begin{abstract}
		The Trust Region Subproblem is a fundamental optimization problem that takes a pivotal role in Trust Region Methods. However, the problem, and variants of it, also arise in quite a few other applications. 
In this article, we present a family of iterative Riemannian optimization algorithms for a variant of the Trust Region Subproblem that replaces the inequality constraint with an equality constraint, and converge to a global optimum.
Our approach uses either a trivial or a non-trivial Riemannian geometry of the search-space, and requires only minimal spectral information about the quadratic component of the objective function. We further show how the theory of Riemannian optimization promotes a deeper understanding of the Trust Region Subproblem and its difficulties, e.g., a deep connection between the Trust Region Subproblem and the problem of finding affine eigenvectors, and a new examination of the so-called hard case in light of the condition number of the Riemannian Hessian operator at a global optimum. Finally, we propose to incorporate preconditioning via a careful selection of a variable Riemannian metric, and establish bounds on the asymptotic convergence rate in terms of how well the preconditioner approximates the input matrix.

	\end{abstract}

	\section{Introduction}\label{sec:intro}
		In this paper, we consider the solution of the following problem, which we term
as the Boundary Trust Region Subproblem (\btrs{}):
\begin{equation}
	\min q\left(\x\right)\coloneqq\frac{1}{2}\x^{\T}\matA\x+\bt\x\st\TNorm{\x}=1, \label{eq:btrs.def}
\end{equation}
where $\matA \in \RR^{n \times n}$ is symmetric and $\b \in \RR^n$.
\btrs{} is closely related to the well known Trust Region Subproblem (\trs{}), which arises in Trust Region Methods:
\begin{equation}
\min q\left(\x\right)\st\TNorm{\x}\leq 1. \label{eq:trs.def}
\end{equation}
Indeed, \btrs{} (\cref{eq:btrs.def}) simply replaces the inequality constraints $\TNorm{\x}\leq 1$ with the equality constraints $\TNorm{\x}=1$. Clearly, these problems coincide whenever a solution for \trs{} can be found on the boundary. Furthermore, an algorithm for solving \btrs{} can be used as a component in an algorithm to solve \trs{}, e.g., in~\cite{Adachi2017}. Alternatively, a simple augmentation trick can be used to translate an $ n $ dimensional \trs{} to an equivalent $ n+1 $ dimensional \btrs{}~\cite{Phan2020}.

While the solution of \trs{} was the initial motivation for our study, our study is also well motivated by the fact that \btrs{} arises in quite a few applications. Indeed, \btrs{} is a form of a constrained eigenvalue problem~\cite{GGM89} that, in turn, arises in machine learning applications such as transductive learning~\cite{Joachims03}, semi-supervised support vector machines~\cite{CSK06}, etc. It also arises when solving quadratically constrained least squares problems~\cite{GM91}, which are closely related to ridge regression. Phan et al. discussed applications of \btrs{} in the context of constrained linear regression and tensor decomposition~\cite{Phan2020}. Finally, we mention recent work on robust function estimation with applications to phase unwrapping~\cite{CT18}.

In many applications, there is a need to solve large-scale instances of \btrs{} or \trs{}, so even if the matrix $\matA$ is accessible and stored in-memory, direct operations such as matrix factorizations are not realistic in terms of running times and/or memory requirements. For this reason, methods that rely solely on matrix-vector products, i.e matrix-free iterative algorithms, are of great interest and appeal when approaching these problems. In this paper we focus on developing matrix-free iterative algorithms for \btrs{} (and \trs{}). We also propose a family of preconditioned iterative algorithms for solving \btrs{} and \trs{}. 

The proposed algorithms are based on  \emph{Riemannian optimization}~\cite{Absil2008}, that is, constrained optimization algorithms that utilize smooth manifold structures of constraint sets. 
Although \btrs{} is a non-convex problem which can have non-global local minimizers, we show that it is possible to find a global solution  using an almost trivial modification of standard Riemannian Gradient Descent\footnote{Standard in the sense that it uses the most natural choice of retraction and Riemannian metric.} without {\em any} spectral information about the matrix $\matA$ (\Cref{sec:solving.naive}). Next, we show that this can be taken one step further, and find a global solution of a \btrs{} with first-order Riemannian algorithms and/or other choices of retraction and/or non-standard Riemannian metric, as long as we have access to the eigenvectors of $\matA$ associated with its smallest eigenvalue (\Cref{sec:solving.prec}); a requirement that is computationally feasible.

It is well known that matrix-free iterative methods may suffer from slow convergence rates, and  that preconditioning can be effective in improving convergence rates of iterative solvers.
Using Riemannian optimization, we are able to perform \emph{Riemannian preconditioning}~\cite{MS16, Shustin2019} by choosing a non-standard~\cite[Eq. 2.2]{edelman1998geometry} metric. Indeed, Riemannian preconditioning introduces a preconditioner by changing the Riemannian metric. We show how to precondition our proposed Riemannian algorithms using an easy-to-factorize approximation $\matM$ of $\matA$. To justify the use of the preconditioner, we present a theoretical analysis that bounds the condition number of the Hessian at the optimum (a useful proxy for assessing convergence rate of Riemannian solvers) in terms of how well $\matM$ approximates $\matA$ (\Cref{subsec:precond-scheme}).

As with any preprocessing, the construction of a preconditioner is expected to have additional computational costs. Our theoretical results are supported by numerical illustrations showing that our preconditioning scheme introduces a speedup large enough in order to result in overall computational costs that are reduced in comparison with iterative schemes based on the standard geometry (\Cref{sec:experiments}).

\subsection{Contributions and organization}
Our work is the first to tackle \btrs{} directly using Riemannian optimization, without re-formulating the problem (see~\Cref{subsec:Preliminaries.related}). 
The main contributions established by viewing this problem from the lens of Riemannian optimization are:
\begin{itemize}
    \item Theoretically, we analyze the possible critical points and their stability, and we explore connections between \trs{} and finding affine eigenvectors using Riemannian optimization theory (\Cref{sec:solving.optimality_and_stat}). 
     In addition, we analyze the easy and hard cases (\Cref{sec:solving.naive}), and show a theoretical relation between the hard case of \btrs{} and the condition number of the Riemannian Hessian at the optimum (\Cref{subsec:precond-scheme}).
    \item Algorithmically, we propose to find a global solution of a \btrs{} via a Riemannian optimization algorithm (\Cref{sec:solving.naive}). 
    Furthermore, we utilize the technique of Riemannian preconditioning~\cite{MS16}, and incorporate preconditioning using a variable Riemannian metric in order to improve the convergence rates of our Riemannian optimization \btrs{} solver (\Cref{subsec:precond-scheme}). Similarly to ~\cite{Shustin2019}, we analyze the effect of preconditioning on the asymptotic convergence by establishing bounds on the condition number of the Riemannian Hessian at the optimum. However, unlike in ~\cite{Shustin2019}, we propose a variable Riemannian metric which adapts itself as iterations progress. Moreover, we design our preconditioner using minimal spectral information about the quadratic term of the objective function, $\matA$, and using matrix sketching techniques which provide efficient computational costs per iteration compared with the use of the exact matrices. In~\Cref{sec:experiments}, 
    we demonstrate the improvement obtained using our Riemannian preconditioning scheme in comparison with naive Riemannian optimization methods without preconditioning.
\end{itemize}

From here on, our text is organized as follows: \Cref{sec:Preliminaries} contains the related work and preliminaries on Riemannian optimization and preconditioning, in \Cref{sec:solving.optimality_and_stat} we study the stationary points of \btrs{} in the Riemannian optimization framework, in \Cref{sec:solving.naive} we propose an adaptation of Riemannian gradient descent which solves \btrs{} globally, and is suitable both for the easy and hard case of \btrs{}, in \Cref{sec:solving.prec} we widen the class of Riemannian solvers which solve globally \btrs{} and utilize it in \Cref{subsec:precond-scheme} to construct and analyze a specific preconditioning scheme, in \Cref{sec:trs} we show how to utilize our solution for \btrs{} for achieving a solution for \trs{}, finally in \Cref{sec:experiments} we illustrate our algorithms for the easy, "almost hard" and hard cases and demonstrate the effect of preconditioning empirically.

	\section{Preliminaries}\label{sec:Preliminaries}
		\subsection{Notation}\label{subsec:Preliminaries.notations}
        
    We denote scalars by lower case Greek letters without subscripts or
    using $x,y,z\dots$. Vectors in $\mathbb{R}^{n}$ are denoted by bold
    lowercase English letters, e.g., $\x,\y,\z\dots$ and matrices by $\matA,\matB,\matC \dots$. The $n \times n$ identity matrix will be denoted by $\matI_{n} $ while the subscript is omitted in cases where the dimension is clear from context.

    Let $\matA \in \RR^{n \times n}$  be a symmetric matrix. We denote its eigenvalues by $\lieig{1}{\matA} \leq \lieig{2}{\matA} \leq \dots \leq \lieig{n}{\matA}$ or simply $\lambda_i \leq \lambda_{i+1}$ where the matrix is clear from the context. We also use $\lmin$ and $\lambda_{\max}$ to denote the minimal and maximal eigenvalue. For any matrix $\matB \in \RR^{n \times n}$, the condition number of $\matB$, denoted by $ \kappa(\matB) $, is defined as the ratio between the largest and smallest singular values of $ \matB $. We say that $\matA \in \RR^{n \times n}$ is symmetric positive definite matrix (\spdm) if $\matA$ is symmetric and all its eigenvalues are strictly positive. In particular, for an \spdm{} matrix, the condition number becomes the ratio between the largest and smallest eigenvalue.

    We denote the $n-1$ dimensional sphere in $\mathbb{R}^{n}$ by $\sph \coloneqq \{ \x \in \RR^{n} ~|~ \x^{\T}\x=1\}$. Recall that $\sph$ is a $ n-1 $ dimensional submanifold of $\RR^n$. Given a smooth function $f : \sph \to \RR$, we use the notation $\bar{f}$ to present \emph{some} smooth extension of $f$ to the entire ambient space of $\RR^n$, that is, $\bar{f}$ refers to \emph{any} member of the equivalence class of smooth functions $g: \RR^n \to \RR$ such that $g(\x) = f(\x)$ for all $\x \in \sph $. We denote the unit norm ball (with respect to the Euclidean norm) in $\mathbb{R}^{n}$ by $\mathbb{B}^{n}$.

    Notation within the domain of optimization algorithms on Riemannian manifolds are consistent with the ones in~\cite{Absil2008}, e.g., general manifolds are denoted using calligraphic uppercase English letters $\mathcal{M},\mathcal{N},\dots$. Similarly, for $\x \in \mathcal{M}$,  the tangent space to $\mathcal{M}$ at $\x$ is denoted by $\tspM{\x}{\mathcal{M}}$, and tangent vectors are denoted by lowercase greek letters with a subscript referring the point for which they correspond, e.g., $\eta_{\x}\in\tspM{\x}{\mathcal{M}}$.

\subsection{Related work}\label{subsec:Preliminaries.related}
    Due to its pivotal role in Trust Region Methods, there has been extensive work on solving \trs{}. 
    There is a variety of classical algorithms to approximate \trs{} solution such as the Cauchy-point algorithm, the dogleg method, two-dimensional subspace minimization and Steihaug's CG algorithm (see \cite[Chapter 4.1]{Nocedal2006} and citations therein). 
    Worth nothing is the seminal work of Mor\'{e} and Sorensen \cite{More1983}, relating solutions of \trs{} to roots of a secular equation. 
    Another classical algorithm for solving \trs{} at large-scale is based on the Lanczos method~\cite{GouldEtAl99}. 

    Recent work by Carmon and Duchi on \trs{} include an analysis of the convergence rate of the Lanczos method ~\cite{Carmon2018}. They also prove lower bounds on the computational cost for any deterministic iterative method, accessing $\matA$ only through matrix-vector products and for which each iteration involves only a single product (i.e., matrix-free algorithms). 
    Beck and Vaisbourd proposed to find  global solutions of \trs{}s by means of first order conic methods (\focm{})~\cite{Beck2018}, and formulated sufficient conditions for such schemes to converge to the global \trs{} minimizer both in easy and hard cases.

    Most previous work on \btrs{} was motivated by \trs{}. Such works addresses \btrs{}  only for the special case where the solutions for both problems coincide. However, there are a few exceptions. Mart{\'{i}}nez  characterized local minimizers of \btrs{} and \trs{}~\cite{Martinez1994} and proposed an algorithm for finding a local but non-global minimizer for \btrs{}, when it exists. 
    This characterization is of  importance when discussing first-order iterative methods as they usually guarantee convergence to a stationary point without any way to distinguish global from local non-global minimizers.
    Given a stationary point of \btrs{} other than its global solution, Lucidi et al. presented a transformation (mapping from a vector to another) for finding a point on the sphere for which the objective value is lower~\cite{Lucidi1998}. If we use a descent algorithm,  Lucidi et al.'s transformation allows us to continue the iteration after converging to a local non-global minimizer (or to any other stationary point). Hager presented an algorithm for solving \btrs{} using a method based on a combination of Newton and Lanczos iterations for Krylov subspace minimization~\cite{Hager2003}. 
    Adachi et al. proposed to solve \btrs{} by solving a  Generalized Eigenvalue Problem of dimension $2n \times 2n$ ~\cite{Adachi2017}. Phan et al. proposed an algorithm for solving \btrs{}, however their algorithm requires a full eigendecomposition of $\matA$~\cite{Phan2020}.

    Our proposed algorithm differs from the aforementioned works in three fundamental ways: 
    1) We consider the use of Riemannian optimization for \trs{}. The only previous work that considered Riemannian optimization for \btrs{} is a recent work by Boumal et al.~\cite{BVB18}, which considers a semidefinite program relaxation which is followed by a Burer-Monteiro relaxation. Their reformulated problem is an optimization problem constrained on two spheres, an $n \times p$ dimensional sphere and a $p$ dimensional sphere, where $p\geq 2$ is a rank parameter~\cite[Section 5.2]{BVB18}. Unlike~\cite{BVB18}, we solve \btrs{} directly via Riemannian optimization. 
    2) Consequently, for \trs{}, our algorithm seeks the solution of an equivalent $ n+1 $ dimensional \btrs{}, from which it is trivial to extract the solution for the original \trs{}.
    3) We incorporate a preconditioner through the approach of Riemannian preconditioning~\cite{MS16}, and not via change-of-variables or preconditioning the solution of linear systems encountered during the optimization (e.g.,~\cite{vandereycken2010riemannian}). Unlike~\cite{MS16}, we motivate the design of our preconditioner via the condition number of the Riemannian Hessian at the optimum.

\subsection{Riemannian Optimization}\label{subsec:Preliminaries.ropt.shpere}
Our proposed algorithms use Riemannian optimization for finding the optimal solution of \btrs{}.
The framework of Riemannian optimization naturally arises when solving optimization problems in which the search space is a smooth manifold~\cite[Chapter 3.1]{Absil2008}. In this section we recall some basic definitions of Riemannian optimization, and establish corresponding
notation. The definitions and notation here are consistent with the ones in~\cite{Absil2008}.

A smooth Riemannian manifold is a differentiable manifold $ \calM $, equipped with a smoothly varying inner product $\gsub{\x} $ operating on the manifold's tangent bundle $ T \calM $, i.e., for any $ \x \in \calM $ the function $ \gsub{\x} : \tspM{\x}{\calM} \times \tspM{\x}{\calM} \rightarrow \RR  $ is a bilinear function on the tangent space to the manifold $ \calM $ at point $ \x $. In turn, this inner product endows a metric function over the tangent space at each point. This inner product is termed the {\em Riemannian metric}.

Riemannian optimization algorithms are derived by generalizing various algorithmic components used in non-Riemannian optimization, and as such are naturally defined on $\RR^n$, to the case of optimization on Riemannian manifolds. For example, a \emph{retraction} \cite[Section 4.1]{Absil2008}, which is a map $R_{\x}:T_{\x}{\cal M}\to{\cal M}$, allows Riemannian optimization algorithms to take a step at point
$\x\in{\cal M}$ in a direction $\xi_{\x}\in T_{\x}{\cal M}$. Two mathematical objects that are important for our discussion are the \emph{Riemannian gradient} and the \emph{Riemannian Hessian}
\cite[Section 3.6 and 5.5]{Absil2008}. 

Once these various components are generalized, many optimization algorithms for smooth
problems are naturally generalized as well. In \cite{Absil2008}, Riemannian gradient,
line-search, Newton method, trust region, and conjugate gradient (CG) methods are
presented. An important example is Riemannian gradient descent, which is given by the following formula:
\begin{equation}
\xkp = R_{\xk}(- t^{(k)} \grad f (\xk ) ), \label{eq:rgd.step.def}
\end{equation}
where $t^{(k)}$ denotes the $k$'th step size. In the above, $\grad f (\xk )$ is the Riemannian gradient of $f$ at $\xk$.
When the step size is chosen via \emph{Armijo's backtracking procedure}, it is guaranteed that all the accumulation points of a sequence generated by Riemannian Gradient Descent are stationary points of $f$ on $\calM$ (vanishing points of the Riemannian gradient), provided $f$ is at least continuously differentiable~\cite[Theorem 4.3.1]{Absil2008}. In general, henceforth, when we discuss  Riemannian Gradient Descent we assume that step sizes are chosen so as to assure that all accumulation points are stationary points (e.g., using Armijo's backtracking procedure).

\subsection{Riemannian Preconditioning on the Sphere}\label{subsec:Preliminaries.rprec.shpere}
The natural way to define a metric on the sphere $\sph$ is by using the standard inner product of its ambient space $\RR^n$: $\bar{g}_{\x}(\eta_\x, \xi_\x) \coloneqq \eta^{\T}_\x \xi_\x$. The sphere $\sph$, as a Riemannian submanifold of $\RR^n$, then inherits the metric in a natural way. With this metric, we have $\gsub{\x}(\eta_\x, \xi_\x) \coloneqq \eta^{\T}_\x \xi_\x$  where $\eta_\x, \xi_\x \in \tsp{\x}$ are given in ambient coordinates.
 
However, Shustin and Avron noticed that in some cases this particular choice of metric may lead to suboptimal performance of iterative algorithms~\cite{Shustin2019}.  For example, when minimizing the Rayleigh quotient defined by an \spdm{} $\matA$, the metric defined by $\matA$ , i.e., $g_{\x}(\eta_\x,\xi_\x)=\eta^{\T}_\x\matA\xi_\x$, was shown to be advantageous~\cite[Section 4]{Shustin2019}.
In general, different problems call for the use of metrics for the form $g_{\x}(\eta_\x,\xi_\x)=\eta^{\T}_\x\matM\xi_\x$ with different $\matM$. 
As the usage of this metric in Riemannian optimization algorithms requires the ability to solve linear systems involving $\matM$, one often wants an $\matM$ that is both easy to invert and closely approximates some optimal (but computationally expensive) metric, e.g., for minimizing $\min_{\x\in\sph}\x^{\T}\matA\x$ we want an easy-to-invert $\matM \approx \matA$. 

Defining the metric on $\sph$ via $\matM$ is an instance of so-called Riemannian Preconditioning~\cite{MS16}.
In our preconditioned iterative algorithms for \btrs{}, a preconditioner is incorporated using Riemannian preconditioning, that is, we use Riemannian optimization on $\sph$ with a non-standard metric. 
However, in contrast to the work by Shustin and Avron~\cite{Shustin2019}, where the metric
is defined by a constant preconditioner $\matM$, our algorithm uses a metric that varies on $\sph$, i.e, a function $g_{\x}(\eta_\x,\xi_\x)=\eta^{\T}_\x\matM_\x\xi_\x$, where for each $\x\in\sph$ the matrix $\matM_{\x}$ is an \spdm{}, and as such it defines a valid inner product on the tangent space to $ \sph $ at $ \x $, and the mapping $\x \mapsto \matM_{\x}$ is smooth on the sphere (smoothness is required in order for $(\sph, g)$ to be a Riemannian manifold).

A summary of Riemannian optimization related objects and their expressions in ambient coordinates is given below in~\Cref{tab:formulas}.

\begin{table}[H]
	\caption{Riemannian optimization related ingredients for
		optimizing on $\protect\sph$ with varying metric. Based on~\cite{Shustin2019}.
		Formulas are given in terms of ambient coordinates. Note that $\x\in\sph\mapsto\matM_{\x}$ is a smooth \spdm{}-valued function. The construction of the objects in the \btrs{} column are given in~\cref{subsec:precond-scheme}  }
	\label{tab:formulas}
	\begin{centering}
	    \scriptsize
		\begin{tabular}{>{\raggedright}m{4.5cm}|>{\raggedright}p{6.7cm}|>{\raggedright}p{3.cm}}
			& {\large $\min_{\x\in\sph}f\left(\x\right)$ \par} & {\large \btrs{} \par}
			\tabularnewline
			\hline
		
			{Tangent space to a point $\x\in\sph$} & {$\tsp{\x}=\left\{ \z\in\R^{n}\,:\,\z^{\T}\x=0\right\} $} & 
		\tabularnewline
			
			{Retraction} & {$R_{\x}(\xi_{\x})\coloneqq\frac{\x+\xi_{\x}}{\TNorm{\x+\xi_{\x}}}$} & 
			\tabularnewline
			&
			&
			\tabularnewline
			{Riemannian metric} & {$g_{\x}\left(\eta_\x,\xi_\x\right)\coloneqq\eta^{\T}_\x\matM_{\x}\xi_\x$}& 
			\tabularnewline
			&
			&
			\tabularnewline
			{Orthogonal projector on $T_{\x}\sph$} & {$\proj{\x}\coloneqq\left(\matI_{n}-1/(\x^{\T}\matM_{\x}^{-1}\x)\matM_{\x}^{-1}\x\x^{\T}\right)$} &
			\tabularnewline
			&
			&
			\tabularnewline
			{Vector transport} & {${\cal T}_{\eta_{\x}}(\xi_{\x} )\coloneqq \proj{R_{\x} (\eta_{\x})} (\xi_{\x})$} & 
			 \tabularnewline
			 &
			 &
			 \tabularnewline
			{Riemannian gradient} & $\grad{f(\x)}=\proj{\x}\matM_{\x}^{-1}\nabla\bar{f}\left(\x\right)$ & {$\proj{\x}\matM_{\x}^{-1}(\matA\x+\b)$}
			\tabularnewline
			&
			&
			\tabularnewline
			{Riemannian Hessian at stationary $\bx$}

			{(i.e., $\grad{f(\bx)=0})$} & {$\Hess f(\bx)[\eta_{\bx}]=\proj{\bx}\matM_{\bx}^{-1}\left[\nabla^{2}\bar{f}\left(\bx\right)-\bx^{\T}\nabla\bar{f}\left(\bx\right))\matI_{n}\right]\eta_{\bx}$} & {$\proj{\bx}\matM_{\bx}^{-1}\left[\matA-\mu_{\bx} \matI_{n}\right]\eta_{\bx}$}
			\tabularnewline
		\end{tabular}
	\end{centering}
\end{table}

	\section{Stationarity in the \btrs{} and Riemannian Optimization}\label{sec:solving.optimality_and_stat}
		Our goal in this section is to understand the set of stationary points of \btrs{}, discuss optimality conditions, and understand how this pertains to solving \btrs{} using plain Riemanniann optimization. Some of the results are closely related to similar results for \trs{}~\cite{Beck2018}, but there are subtle differences. 

Recall, that for a Riemannian manifold $ \mathcal{M} $, a stationary point $\x \in \mathcal{M}$ of a smooth scalar function $ f : \mathcal{M} \rightarrow \RR $ is a point for which the Riemannian gradient vanishes ($ \grad f (\x) = 0$)~\cite{Absil2008}. In this section we analyze stationarity of $\x \in \sph$ for $q$. Since vanishing points of the Riemanninan gradient are invariant to the choice of the metric (as locally, the Riemannian metric is an inner product on the tangent space), we can analyze stationarity with {\em any} Riemannian metric of our choice. In this section, we view $\sph$ as a Riemannian submanifold of $\RR^n$ endowed with the dot product as the Riemannian metric. 

The following proposition characterizes the stationary points of \btrs{}. 
The result is already known ~\cite{Martinez1994,Chambers2001,Luenberger1986,More1983}. 
Nevertheless, we provide a new proof, which is based on Riemannian optimization tools. 

\begin{proposition}\label{thrm:btrs.stat}
	A point $ \x \in \sph $ is a stationary point of \btrs{} if and only if there exists $ \mu_\x \in \RR $ such that 
	\begin{equation}\label{eq:btrs.stat.b}
		(\matA - \mu_\x \matI) \x = -\b. 
	\end{equation}
	When such is the case, $\mu_\x$ is unique, and
	\begin{equation}\label{eq:btrs.stat.mu}
	    \mu_\x = \x^{\T}\matA\x+\b^{\T}\x .
	\end{equation}
\end{proposition}

\begin{proof}[Proof of \Cref{thrm:btrs.stat}]
Since we are viewing $\sph$ as a Riemannian submanifold of $\RR^n$ equipped with usual dot product, we have 
$$
\grad{q(\x)} = \proj{\x}\nabla\bar{q}(\x) = (\matI_n - \x \x^\T)(\matA\x + \b) 
= \matA \x + \b - (\x^{\T}\matA\x+\b^{\T}\x)\x ~~.
$$

In the above, $\nabla \bar{q}(\x)$ is the Euclidean gradient of $\bar{q}$ at $\x$, and $\matP_\x \coloneqq \matI_n - \x \x^\T$ is the projection matrix on $\tsp{\x}$ with respect to the Euclidean inner product (the dot product). 
The fact that $\grad{q(\x)} = \proj{\x}\nabla\bar{q}(\x)$ is due to a generic result on the Riemannian gradient of a function on a Riemannian submanifold~\cite{Absil2008}. 
Existence of $\mu_\x$ and the formula given for it (\Cref{eq:btrs.stat.mu}) now follows by equating $\grad{q(\x)}=0$. The converse follows from substituting \Cref{eq:btrs.stat.mu} in \Cref{eq:btrs.stat.b}.

As for uniqueness, if there were two $\mu_{\x,1}$ and $\mu_{\x,2}$ for which \Cref{eq:btrs.stat.b} holds, then obviously $\mu_{\x,1} \x = \mu_{\x,2} \x$. 
Since $\x \in \sph$ we have $\mu_{\x,1} = \mu_{\x,2}$.
\end{proof}

A similar claim holds for \trs{}~\cite{Beck2018}, however for \trs{} we always have $\mu_\x \leq 0$ while for \btrs{} it is possible that $\mu_\x > 0$ (however, this may happen only if $\matA$ is positive definite). The set of pairs $ (\mu_{\x}, \x) $ where  $\x $ is stationary point of \btrs{} is exactly the set of \emph{KKT pairs} for \btrs{}~\cite{Lucidi1998}. It is also the case that any stationary point on $\sph$ of the associated \trs{} is also a stationary point of \btrs{}, but the converse does not always hold. 

In the special case where $\b = 0$, the \btrs{}'s objective function $q(\x)$ is the Rayleigh quotient, and the stationary points are the eigenvectors. In this case, the stationarity conditions reduce to $\matA \x = \mu_\x \x $, so $\mu_\x $ is the corresponding eigenvalue. When $\b \neq 0$, we can still view a stationary $\x $ as an eigenvector, but of an affine transformation instead of a linear one. Indeed, consider the affine transformation $\matT(\v) = \matA\v + \b$. 
We have that $ \matT(\x) = \mu_\x \x $, i.e., $(\mu_\x,\x) $ behaves like an eigenpair of $\matT$.
Furthermore, any affine transformation $ \matS : \RR^n \rightarrow \RR^n $ can be written as $ \matS(\v) = \matA_{\matS} \v + \b_{\matS} $ for some $ (\matA_{\matS}, \b_{\matS}) $. This motivates the following definition which echoes previous observations, e.g. see~\cite{GGM89}, regarding the stationary points of \btrs{}. 

\begin{definition}[Affine Eigenpairs]\label{def:affine.pair}
	Let $\matA \in \RR^{n \times n}$  and  $\b \in \RR^n$. 
	We say that $\mu \in \RR$ is an \emph{affine eigenvalue} of $(\matA,\b)$ if there exists an 
	$\x \in \sph$ such that $\matA\x +\b = \mu \x$. Such $\x\in \sph$ is the \emph{affine eigenvector} associated with $\mu$.
	We call the pair $(\mu,\x)$ an \emph{affine eigenpair}.
\end{definition}

Let us define $\mu_\x \coloneqq \x^{\T}\matA\x+\b^{\T}\x$ for any $\x $ (not just stationary $\x $). If $\b = 0$, for any $\x \in \sph $ the quantity $\mu_\x$ is a Rayleigh quotient of $\matA$. Since $\mu_\x $ plays a similar role for the affine eigenvalues as the Rayleigh quotient plays for the (regular) eigenvalues, we refer to $\mu_\x $ the {\em affine Rayleigh quotient} of $\x $ with respect to $(\matA, \b)$. Like the standard Rayleigh quotient, $\mu_\x$ provides the "best guess"
for the affine eigenvalue, given an approximate affine eigenvector $\x$ since $\mu_\x = \arg\min_\mu \TNorm{\matA \x + \b -\mu \x}$.

\begin{corollary}\label{lem:stat-affine-eval}
	A point $ \x \in \sph $ is a stationary \btrs{} point if and only if $ \x $ is an affine eigenvector of $ (\matA, \b) $, and its associated affine eigenvalue is the affine Rayleigh quotient $ \mu_{\x}$.
\end{corollary}


In their work, Mor\'{e} and Sorensen have shown that when $\b\neq 0$ the affine eigenvalues of $ (\matA, \b) $  are the roots of a secular equation~\cite{More1983}\footnote{We caution that~\cite{More1983} does not use the term "affine eigenvalues".}. As Lucidi et al. later noted, this implies that there is at most one affine eigenvalue smaller or equal to the minimal eigenvalue $ \lamsub{1} $ of $\matA$,  at most two affine eigenvalues between each two distinct eigenvalues, and exactly one affine eigenvalue larger or equal to the largest eigenvalue $ \lamsub{n} $ of $\matA$~\cite{Lucidi1998}. In a sense, when $\matA$ is symmetric there
are at most two affine eigenvalues per each regular eigenvalue, one smaller than it and one larger than it. When $\b=0$, these two affine eigenvalues coincide and have two different eigenvectors, which are the reflection of each other. When $\b$ is perturbed, the affine eigenvectors bifurcate, and when $\b $ is large enough there might fail to be a root of the secular equation (and the affine eigenvalue disappears). 

The following lemma is a compilation of multiple results from \cite{Martinez1994, Lucidi1998} and relates affine eigenpairs to local and global \btrs{} minimizers. Right afterwards we state and prove a refinement of the first clause of the lemma. 

\begin{lemma}[Combining multiple results from~\cite{Martinez1994, Lucidi1998}]~\label[lem]{lemma:affine_pairs_and_minimizers} 
	Let $ \matA \in \RR^{n \times n}$ be a symmetric matrix and $ \b \in \RR^n $. The following statements hold: 
	\begin{lemlist}
		\item Any global \btrs{} solution is an affine eigenvector $ \xopt $ associated with the smallest affine eigenvalue $ \muopt$ and vice versa. We also have $\muopt \leq \lamsub{1}(\matA)$. \label{lemma:globmin_affine_pairs_and_minimizers}
		\item Any stationary point which is not a global solution is an affine eigenvector $ \bx $ associated with an affine eigenvalue $\bar{\mu}$ for which $\bar{\mu} >  \lamsub{1}(\matA)$\label{lemma:nonglobmin.affine.pairs.minimizers} \footnote{This statement is simple corollary of Lemma 2.2 in~\cite{Martinez1994}.}.
		\item If the minimal eigenvalue of $ \matA $ is simple, there are at most two distinct local \btrs{} minimizers~\cite{Lucidi1998}, otherwise, any local \btrs{} minimizer is a global solution~\cite[Lemma 2.2]{Martinez1994}.
		\item In case $ \b \perp \u $ for some $ \u \in \sph $ such that $\matA \u = \lamsub{1}(\matA) \u$, any local \btrs{} minimizer is a global one. \label{lemma:hardcase_localisglobal_affine_pairs_and_minimizers}
		\item Let $ \bar{\x} \in \sph$ be a local but non-global \btrs{} minimizer,  then the affine eigenvalue $ \bar{\mu} $ associated with $ \bar{\x} $ is the second smallest affine eigenvalue, with $ \lamsub{1}(\matA) <  \bar{\mu} < \lamsub{2}(\matA) $. 
	\end{lemlist}
\end{lemma}

\begin{proposition}\label{prop:affine-strict}
	Suppose that $ \ut \b \neq 0 $ for some $ \u \in \sph $ such that $\matA \u = \lamsub{1}(\matA) \u$. 
	Let  $ \muopt$ be the smallest affine eigenvalue.
	Then, $ \muopt < \lamsub{1}(\matA)$. 
\end{proposition}

\begin{proof}
	Let $\xopt$ be an affine eigenvector associated with  $ \muopt$ (i.e., a global minimizer). By~\cref{eq:btrs.stat.b}, for any eigenvector $ \v $ of $ \matA $ such that $ \matA \v = \lambda \v $, it holds that 
	\begin{align*}
	\vt \b &= -\vt (\matA - \muopt \matI_n)\xopt  \\
	&= -(\lambda - \muopt)\vt \xopt\ , 
	\end{align*}
	and in particular we have that $ (\lamsub{1}(\matA) - \muopt) \ut \xopt = - \ut \b \neq 0 $,
	so $ \lamsub{1}(\matA) \neq \muopt $. 
	Since we already know from \Cref{lemma:globmin_affine_pairs_and_minimizers} that $ \muopt \leq \lamsub{1}(\matA) $, we conclude that $ \muopt < \lamsub{1}(\matA) $.
\end{proof}

In general, we can expect a first order Riemannian optimization method to converge to a stationary point, as it is the case with Riemannian Gradient Descent.
The upshot of~\cref{lemma:affine_pairs_and_minimizers} is that we want it to converge to a stationary point whose corresponding affine eigenvalue $\mu_{\x}$ is small, and in particular we want it to converge to the vector associated with the smallest affine Rayleigh quotient.
A key property of first order optimization methods in general, is that given reasonable initialization point and choice of step size the iterations will converge to a stable stationary point; see~\cite[Theorem 4.3.1]{Absil2008} and~\cite[Chapter 4]{boumal2022intromanifolds}. This motivates a study of which affine eigenvectors are stable stationary points.


\if01
\begin{definition}[Riemannian Gradient Descent]\label{def:RGD}
	Let $ \mathcal{M}  $ be a Riemannian manifold and $ f:\mathcal{M} \rightarrow \RR $ a smooth function.
	Given a retraction $ R: \tspM{{}}{\mathcal{M}} \rightarrow \mathcal{M}$ and a starting point $ \x_0 \in \mathcal{M} $ the Riemannian Gradient Descent algorithm forms a sequence of iterates defined by the following step
	\begin{equation}
		\x_{k+1} \gets R_{\x_k} (- \alpha_k \grad f (\x_k))\label{eq:RGD.step}
	\end{equation}
	Where $ \alpha_k $ are the sizes of steps taken at each iterate, and $ - \alpha_k \grad f (\x_k) $ considered as vectors in $ \tspM{\x_k}{\mathcal{M}} $.
	See~\cite[Algorithm 4.1 and definitions therein]{boumal2020intromanifolds} for detailed description.
\end{definition}
\begin{rem}
	In this work, we consider only the "canonical" retraction to the sphere that is length-one normalization.
	With this particular choice of retraction, it is possible to express the retraction at any point using extrinsic coordinates of the ambient space by
	\begin{equation}
		R_{\x} (\eta_\x) \coloneqq (\x + \eta_\x) / \TNorm{\x + \eta_\x} \qquad \forall \x \in \sph, \eta_\x \in \tspM{\x}{\sph} \label{eq:canonical.retraction.ext}
	\end{equation}
	Note that $ \eta_\x $ is interpreted as a "deviation from origin" of the vector space $ \tspM{\x}{\sph} $, which is translated to "deviation from $ \x $"  in $ \RR^n $ (See~\cite[Chapters 3.5 and 4.1]{Absil2008} for a thorough discussion on the subject)
\end{rem}

	The following is a slight strengthening of clause (i).
	\begin{proposition}
		Any stationary point, which is not a global solution is an affine eigenvector $ \x $ associated with with an affine eigenvalue $\mu$ for which $\mu >  \lamsub{1}(\matA)$.
	\end{proposition}
	\begin{proof}
		We divide the proof into two cases: there exists a $\u$ such that $\matA\u = \lamsub{1}\b$ and $\ut \b \neq 0$ and when this does not hold.
		
		 Consider the first case. We show that $\lamsub{1}$ cannot be an affine eigenvalue, so obviously since for every stationary $\x$ that is not a global solution we have $\mu_\x \geq \lamsub{1}$ then $\mu_x$, which must be an affine eigenvalue, must have $\mu_\x > \lamsub{1}$. To see this, assume by contradiction that $\y$ is an affine eigenvector associated with the the affine eigenvalue $lamsub{1}$. Then,
		 $$
		 (A - \lamsub{1}\matI_n)\y = -\b.
		 $$ 
		 Pre-multiply by $\ut$ to get $\ut b = 0$, which is a contradiction. 
		 
		 Uria complete 
	\end{proof}
	We see from the last lemma that there exists a dividing line.... It is customary in the literature on solving \trs{} to refer the former as the {\em easy case} and the latter as the {\em hard case}. We adopt a similar terminology for \btrs{}. 
\fi

The following theorem classifies the stationary points of \btrs{} according to their stability or instability with respect to Riemannian Gradient Descent~\cite[Algorithm 4.1]{boumal2022intromanifolds}. We follow the definitions of~\cite[Section 4.4]{Absil2008} for \emph{stable}, \emph{asymptotically stable}, and \emph{unstable} fixed points. In other words, fixed points for which iterations in a neighborhood of it stay in some neighborhood, converge to the fixed point, or leave the neighborhood correspondingly.
This result helps us understand how plain Riemannian optimization for \btrs{} behaves, and  which among the stationary points of \btrs{} are  unstable, thus reducing the number of probable outcomes of the algorithm.
Although it formally applies only to a specific algorithm, we believe it is indicative for the behavior of other Riemannian first order methods (e.g., Riemannian CG).

\begin{thm}\label{lemma:uria.3.globmin.attractor}
	Let $\{\x_{k}\}$, be an infinite sequence of iterates generated by 
	Riemanninan Gradient Descent as described in~\cite[Algorithm 4.1]{boumal2022intromanifolds} on $q(\x)$. Then the following holds:
	\begin{thmlist}
		\item Every accumulation point of $\{\x_{k}\}$ is an affine eigenvector of $(\matA,\b)$. \label{lemma:attractors.item.accumulation}
		\item The set of affine eigenvectors associated with the minimal affine eigenvalue $ \muopt $ is comprised of stable fixed points of that iteration.\label{lemma:attractors.item.stableglobmin}
		\item In the case $ \muopt < \lamsub{1}(\matA) $, then the affine eigenvector associated with $ \muopt $ is unique, and is an asymptotically stable fixed point. In particular, this occurs when there exists an eigenvector $ \u $ of $ \matA $ such that $ \matA \u = \lamsub{1}(\matA) \u $ for which $ \ut \b \neq 0 $.\label{lemma:attractors.item.assympstableglobmin}
		\item Any affine eigenvector associated with an affine eigenvalue $ \bar{\mu} $ greater than the second smallest affine eigenvalue $ \mu_2 $, i.e., $\bar{\mu} > \mu_2$, is an unstable fixed point.\label{lemma:attractors.item.unstable}
	\end{thmlist}
	
\end{thm}

Before proving the theorem, we first prove a couple of auxiliary results.



\begin{lemma}[\texorpdfstring{Expansion of~\cite[Lemma 3.1]{Lucidi1998}}\nobreak]
	\label{lemma:uria.1.lucidi}
	Let $ (\bar{\mu} , \bar{\x}) $ and $ (\hat{\mu}, \hat{\x}) $ be two affine eigenpairs, then $ \bar{\mu} = \hat{\mu} $ if and only if $ q(\bar{\x}) = q(\hat{\x}) $.
\end{lemma}

\begin{proof}
	The fact that $\bar{\mu} = \hat{\mu}$ implies that $ q(\bar{\x}) = q(\hat{\x}) $ is proved in~\cite[Lemma 3.1]{Lucidi1998}.
	
	For the other direction, assume $ q(\hat{\x})  = q(\bar{\x}) $. Notice it is always the case that $2q(\x) = \mu_\x + \bt \x$ so we have 
	\begin{equation}
	\bar{\mu} + \b^\T \bar{\x} = \hat{\mu} + \b^\T \hat{\x}. \label{eq:obj.iff.mu.1} 
	\end{equation}
	Since both $\bar{\x}$ and $\hat{\x}$ are affine eigenpairs, we have $ \b = - (\matA - \hat{\mu} \matI_n) \hat{\x} = - (\matA - \bar{\mu} \matI_n) \bar{\x}  $ and we can re-write \Cref{eq:obj.iff.mu.1} as 
	$$
	\bar{\mu} - \hat{\x}^\T (\matA - \hat{\mu} \matI_n) \bar{\x} = \hat{\mu} -  \bar{\x}^\T (\matA - \bar{\mu} \matI_n) \hat{\x} ,
	$$
	which can be reduced to 
	$$
	\bar{\mu} - \hat{\mu} = (\bar{\mu} - \hat{\mu})  \hat{\x}^\T \bar{\x}.
	$$
	Now, for this equation to hold we either have  $\bar{\mu} = \hat{\mu}$ 
	(in which case we are done) or $ \hat{\x}^\T \bar{\x} = 1 $. For the latter, since both $\hat{\x}$ and $\bar{\x}$ have unit norm, we must have $\hat{\x}=\bar{\x}$ and again we have $\bar{\mu} = \hat{\mu}$ (the affine eigenvalue corresponding to an affine eigenvector is unique).
\end{proof}

\begin{lemma}\label{lemma:uria.2.disjoint.compacts}
	For an affine eigenvalue $\mu$ denote 
	$$ \mathcal{L}_\mu \coloneqq \{ \x \in \sph\,|\, \matA \x + \b = \mu \x \} ,$$
	(i.e., $\mathcal{L}_\mu$ is the set of affine eigenvectors corresponding to $\mu$).
	We have $ \dist(\mathcal{L}_\xi, \mathcal{L}_\nu) = 0 $ if and only if $ \xi = \nu $, where 
	$$ 
	\dist(\mathcal{L}_\xi, \mathcal{L}_\nu) \coloneqq \inf \{ \TNorm{\x_\xi - \x_\nu}~|~\x_\xi \in \mathcal{L}_\xi,  \x_\nu \in \mathcal{L}_\nu \} .
	$$
\end{lemma}

\begin{rem*}
	If $\matA$ is symmetric and $\mu$ is an affine eigenvalue that is not a (standard) eigenvalue, then it is easy to show that the corresponding affine eigenvector is unique, and $\mathcal{L}_\mu$ contains a single point. However, if the affine eigenvalue $\mu$ is also an eigenvalue, and that eigenvalue is not simple, then the set $\mathcal{L}_\mu$ is not single point. 
\end{rem*}

\begin{proof}[Proof of \Cref{lemma:uria.2.disjoint.compacts}]
	Suppose $ \dist(\mathcal{L}_\xi, \mathcal{L}_\nu) = 0 $. So, without loss of generality, there exist a sequence $ \{\x_\xi^i \}_{i=1}^{\infty} $ of points in $\mathcal{L}_\xi$ such that $ \lim_{i \rightarrow \infty} \x_\xi^i = \x_\nu \in \mathcal{L}_\nu $ (where we used the fact that $\mathcal{L}_\nu$ is closed). Since $q$ is continuous we find that $\lim_{i \rightarrow \infty} q(\x_\xi^i) = q(\x_\nu)$. However, Lemma~\ref{lemma:uria.1.lucidi} implies that $q(\x_\xi^i)$ is constant for all $i$ since all $\x_\xi^i$s are affine eigenvectors of the same affine eigenvalue, which in turn implies that $\lim_{i \rightarrow \infty} q(\x_\xi^i) = q(\x_\xi^1)$. We found that $q(\x_\xi^1)=q(\x_\nu)$ and \Cref{lemma:uria.1.lucidi} now implies that $\xi=\nu$. 
\end{proof}


\begin{proof}[Proof of~\cref{lemma:uria.3.globmin.attractor}]
	\Cref{lemma:attractors.item.accumulation} follows directly from the convergence analysis of Riemannian Gradient Descent (~\cite[Theorem 4.3.1]{Absil2008},~\cite[Propositions 4.7, Corollary 4.9, and Corollary 4.13]{boumal2022intromanifolds}), and~\Cref{lem:stat-affine-eval}.
	
	For~\Cref{lemma:attractors.item.stableglobmin}, we show that any neighborhood $ \mathcal{U} \subset \sph $ containing the set of affine eigenvectors associated with the minimal affine eigenvalue, there exists a non-empty level-set contained in $ \mathcal{U} $ in which the only stationary points are affine eigenvectors corresponding to the minimal affine eigenvalue.
	
	Let $ \mathcal{L}_{\star} $ denote the set of affine eigenvectors associated with the minimal affine eigenvalue $ \muopt $, and let $ \Lambda $ be the set of affine eigenvalues $ \mu  $ with $ \mu > \muopt $.
	By~\cite[Proposition 3.2]{Lucidi1998} the set $ \Lambda $ is finite, and thus 
	\[
	\mathcal{L}_{\Lambda} \coloneqq \bigcup_{\mu \in \Lambda} \mathcal{L}_\mu ,
	\]
	is compact, since it is a finite union of compact sets. 
	For any neighborhood $ \mathcal{U} \subset \sph$ containing $ \mathcal{L}_\star $, write 
	\[
	l_1 \coloneqq \inf_{\x \in \sph \smallsetminus \mathcal{U} } q(\x) ,
	\]
	Note that for any $ l < l_1 $, the level set of points $ \x $ for which $ q(\x) \leq l $ is a subset of $ \mathcal{U} $.
	
	By~\Cref{lemma:globmin_affine_pairs_and_minimizers}, it holds that $ q(\v) = \min_{\x \in \sph} q(\x) $ for all $ \v \in \mathcal{L}_\star $. 
	Write $ q_\star \coloneqq \min_{\x \in \sph} q(\x) $.
	Now let $ l_2 \coloneqq \min_{\x \in \mathcal{L}_\Lambda} q(\x)$, and note that for any $ l < l_2 $ the intersection of $ \mathcal{L}_\Lambda $ and the level set of points $ \x $ such that $ q(\x) \leq l $ is empty. 
	Define $ l \coloneqq \min \left\{  l_1 + q_\star  ,  l_2 + q_\star   \right\} / 2 $, and observe that $ q_\star < l < l_i   $ for both $ i=1,2 $, and let
	$$
	\mathcal{L} \coloneqq \left\{ \x \in \sph ~|~ q(\x) \leq l  \right\} .
	$$
	By construction, we have that $ \mathcal{L} \subseteq \mathcal{U} \smallsetminus \mathcal{L}_\Lambda $. 
	So we showed that any neighborhood $ \mathcal{U} $ of $ \mathcal{L}_\star $ contains a sub-level set $ \mathcal{L} \supseteq  \mathcal{L}_\star$ such that $ \x \in \mathcal{L} $ is a stationary point if and only if $ (\matA - \muopt \matI_{n}) \x = -\b $, which concludes~\Cref{lemma:attractors.item.stableglobmin}.
	
	If in addition  $ \muopt < \lamsub{1}(\matA) $ we have that the global minimizer is unique, in which case any descent mapping starting at $ \x_0 \in \mathcal{L} $ will surely converge to the only critical point in that level-set, which is the affine eigenvector $ \xopt $ associated with $ \muopt $, thus~\Cref{lemma:attractors.item.assympstableglobmin} holds.
	Note that in case $ \bt \u \neq 0 $ for some $ \u $ such that $ \matA \u = \lamsub{1}(\matA) \u $, it is clear that $ \muopt < \lamsub{1}(\matA) $ by~\Cref{prop:affine-strict}.
	
	As for~\Cref{lemma:attractors.item.unstable}, we know that in addition to the global minimizer, there is (potentially) only one more local minimizer that is not global, which, if exists, is an affine eigenvector associated with the second smallest affine eigenvalue.
	As affine eigenvectors corresponding to values $ \bar{\mu} > \mu_2 $ cannot be a local minimizer, 
	and \Cref{lemma:uria.2.disjoint.compacts} ensures that every such affine eigenvector has a compact neighborhood where every other stationary point in that neighborhood has the same objective value, then according to~\cite[Theorem 4.4.1]{Absil2008} this affine eigenvector must be an unstable fixed points.
\end{proof}


%
Thus, for most initial points, we can expect first order Riemannian optimization methods to converge to one of at most two local minimizers. One of the local minimizers is the global minimizer, but the other one might not be. The local non-global minimizer corresponds to a small affine eigenvalue, and heuristically it should have a not too bad objective value. We see that plain Riemannian optimization is not a bad choice. Nevertheless, we are interested in methods which find a global solution. In subsequent sections we propose Riemannian optimization methods that converge to a global optimum.

	\section{First-Order Riemannian \btrs{} Solver which Converges to a Global Optimum}\label{sec:solving.naive} %
		In this section we present a solver for \btrs{} that uses Riemannian optimization and finds a global solution of a \btrs{}. 
Our proposed algorithm is listed in~\Cref{alg:naive.convergent}.
Remarkably, our algorithm requires no spectral information on the matrix $\matA$.
Similarly to ~\cite{Beck2018}, our analysis identifies sufficient optimality conditions for isolating the global solution for each of the two \btrs{} cases.
Hence the double-start strategy employed by~\Cref{alg:naive.convergent}; without any assumptions regarding the current \btrs{} case, the Riemannian optimization is initiated from two distinct starting points (corresponding to each of the optimality conditions).
Yet, the underlying idea of~\Cref{alg:naive.convergent} differs from those presented in~\cite{Beck2018}: while~\cite{Beck2018} relies on  \focm{} steps (concretely - Projected/Conditional Gradient methods), our proposed algorithm uses Riemannian Gradient Descent.
In fact, the use of Projected Gradient Descent for global solution of \btrs{} on the sphere does require knowledge about $\matA$'s spectral properties.

Our algorithm uses Riemannian Gradient Descent, but fixes a specific Riemannian metric and a particular retraction. 
The Riemannian metric is simply obtained by viewing $\sph$ as a submanifold of $\RR^n$ which, is viewed as an inner product space equipped with the usual dot product. In this context, the Riemannian gradient of the objective $q$ ~\Cref{eq:btrs.def} at a point $\x \in \sph$ is given by
	\begin{equation}
	\grad q (\x) = \matA \x + \b - (\x^{\T} \matA \x ) \x - (\x^{\T} \b) \x. \label{eq:se.naive.rgrad}
	\end{equation}
For the retraction, we project to $\sph$ by scaling:
$$
R_\x(\xi_\x) = \frac{\x + \xi_\x}{\TNorm{\x + \xi_\x}} ~~.
$$

Thus, our algorithm employs the following iteration:
\begin{equation}\label{eq:foro.step}
	\xkp = \frac{\xk - t^{(k)} \grad q (\xk)}{\TNorm{\xk - t^{(k)} \grad q (\xk}}\ ,
\end{equation}
where $t^{(k)}$ is a positive step-size (e.g., chosen by backtracking line search).


As is customary for \trs{}, we can classify instances of \btrs{} into two cases, "easy cases" and the "hard cases", based on the relation between $\b$ and $\matA$. 
The easy case is when there exists an eigenvector $\u$ such that $\matA\u=\lamsub{1}(\matA)\u$  for which $\ut\b\neq0$. 
The hard case is when no such $\u$ exists. 


The strategy employed by the algorithm is similar to the one in~\cite{Beck2018}: we identify two sets, $S_E$ and $S_H$, such that in the easy case (respectively, the hard case) the global solution is the only stationary point in $S_E$ (respectively, $S_H$). We then derive a condition that ensures that if the iterations starts in $S_E$ (respectively, $S_H$), then all accumulation points of the generated sequence of iterations are in $S_E$ (respectively, $S_H$). By executing two iterations, one starting in $S_E$ and another starting in $S_H$, we cover both cases.

We remark that the definition of $S_E$ and $S_H$ is the same as the one given in~\cite{Beck2018} for \trs{}. The proofs in this section are also similar to proofs of analogous claims in~\cite{Beck2018}. However, nontrivial adjustments for \btrs{} were needed.

We begin with the easy case. Let us define $S_E$:
	\begin{equation}
	S_E \coloneqq \{ \x \in \sph~|~ (\vt \b)(\vt \x) \leq 0 \quad \forall \v \st \matA\v = \lamsub{1}(\matA)\v \}. \label{eq:SE.ball.eq}
	\end{equation}
The following lemma shows that the global optimum belongs to $S_E$, and if we are in the easy case, the global optimum is the only stationary point in $S_E$.
\begin{lemma}\label{lemma:se.sphere.cond}
	Let $ \xopt $ be a global minimum of \btrs{}. Then $ \xopt \in  S_E $. Futhermore,  if there exists a $\u$ such that $\matA\u = \lamsub{1}(\matA)\u$ and $\b^\T \u \neq 0$, then $\xopt $ is the only stationary point in $S_E$. 
\end{lemma}

\begin{proof}
First, let us show that $ \xopt \in S_E$. 
 The characterization of stationary points of \btrs{} enables us to write
	\[
	\b = - (\matA - \muopt \matI) \xopt\ ,
	\]
	where $\muopt$ is the affine eigenvalue associated with $\xopt$. Consider an eigenvector $\v$ corresponding to the smallest eigenvalue of $\matA$. Pre-multiplying by $ (\vt \xopt) \vt $ results in 
	\[
	(\vt \xopt) (\vt \b) = - (\lamsub{1} - \muopt) (\vt \xopt)^2.
	\]
	According to \Cref{lemma:globmin_affine_pairs_and_minimizers} we have $ \muopt \leq \lamsub{1} $ and thus $ (\vt \xopt) (\vt \b) \leq 0$. This holds for an arbitrary eigenvector $\v$ corresponding to the minimal eigenvalue so $\xopt \in S_E$.
	
%
	
	Next, we show that in the easy case, $\xopt$ is the only stationary point in $S_E$.
	Let $\bar{\x}$ be a stationary point that is not a global minimum, and let $\bar{\mu}$ be its associated affine eigenvalue. Let $\u$ be an eigenvector corresponding to the smallest eigenvalue for which $\bt \u \neq 0$ (we assumed that such an eigenvector exists). 
	First, we claim that $\ut \bar{\x} \neq 0$. To see this, recall that 
	$$
	(\matA - \bar{\mu}\matI_n)\bar{\x} = -\b,
	$$
	and pre-multiply this equation by $\ut$ to get 
	$$
	(\lamsub{1} - \bar{\mu})(\ut\bar{\x}) = -\ut \b.
	$$
	Since $ \ut \b \neq 0 $, both $\ut\bar{\x}\neq 0$ and  $\lamsub{1} \neq \bar{\mu}$ must hold. Moreover, due to \Cref{lemma:affine_pairs_and_minimizers} we have $\bar{\mu} > \lamsub{1}(\matA)$. Pre-multiply by $\ut\bar{\x}$ to get:
	$$
	- (\ut\bar{\x})(\ut \b) = (\lamsub{1} - \bar{\mu})(\ut\bar{\x})^2.
	$$
	We have $(\ut\bar{\x})^2 > 0$ and $(\lamsub{1} - \bar{\mu}) < 0$ so $- (\ut\bar{\x})(\ut \b) < 0$, violating the definition of $S_E$, thus we have $\bar{\x} \not\in S_E$.
\end{proof}

Next, we show that if the step-size is small enough (smaller than $1/\TNorm{\b}$), then if $\xk \in S_E$ we also have $\xkp\in S_E$. Thus, iterations that start in $S_E$ stay in $S_E$. This is obvious in the hard case (where $S_E = \sph $ ). However, it proves useful in the easy case due to the previous lemma and the fact that $S_E$ is closed. 
\begin{lemma}\label{lemma:unprec.se.step}
	Consider \Cref{eq:foro.step}. Provided that $t^{(k)} < 1/\TNorm{\b}$, if $\xk \in S_E$ then we also have $\xkp \in S_E$.
\end{lemma}
\begin{proof}
		We can write 
	\begin{equation}\label{eq:xkpstep}
	\xkp = \theta_1 \xk + \theta_2 ( - \grad q ( \xk ) ),
	\end{equation}
	where 
	\begin{align*}
	\theta_1 = \TNorm{ \xk - t^{(k)} \grad q (\xk) }^{-1} \quad \theta_2 = t^{(k)} \TNorm{ \xk - t^{(k)} \grad q (\xk) }^{-1}\ .
	\end{align*}
	Now let $\v$ be an eigenvector of $\matA$, corresponding to its smallest eigenvalue. By~\cref{eq:se.naive.rgrad} we have:
	\begin{align*}
	\vt \grad q (\xk) &= \vt \matA \xk + \vt \b - (\xkt \matA \xk) \vt \xk - (\xkt \b) \vt \xk\\
	&= (\lamsub{1} - \xkt \matA \xk)   \vt \xk + \vt \b - (\xkt \b) \vt \xk\ , 
	\end{align*}
	and hence: 
	\begin{align}
	\vt\b ( - \vt \grad q ( \xk ) ) &= \vt\b ( - ( \lamsub{1}  - \xkt \matA \xk ) \vt\xk - \vt\b + ( \xkt \b ) \vt\xk )\nonumber\\
	&\leq   (\vt \b )  (\xkt \b) (\vt \xk),\label{eq:naive.se.ineq}
	\end{align} 
	where we used the fact that $\xkt \matA \xk \geq \lamsub{1}$ and $(\vt \b)(\vt \xk) \leq 0$ since $\xk \in S_E$. Now, further write:
	\begin{align}
	(\vt \b)( \vt \xkp )&= \theta_1(\vt \b)(\vt \xk) + \theta_2 ( \vt \b )( - \vt \grad q(\xk) ) \label{eq:ut_xkp_ineq}\\
	&\leq \theta_1 (\vt \b)(\vt \xk) +  \theta_2  (\vt \b )  (\xk^{\T} \b) (\vt \xk) \nonumber \\
	&= (\theta_1 +  \theta_2  (\xk^{\T} \b) ) (\vt \b )(\vt \xk).  \nonumber 
	\end{align}
	Since $\xk \in S_E$  we have $ (\vt \b )(\vt \xk)  \leq 0 $, we are left with showing that $ \theta_1 +  \theta_2  (\xk^{\T} \b)  \geq 0 $, which obviously occurs when $ 1 + t^{(k)} (\xk^{\T} \b) \geq 0 $.
	Since $ \TNorm{\xk} = 1 $, we have that $ \xk^{\T} \b \geq - \TNorm{\b} $. 
	By assuming that $ t^{(k)} < 1 / \TNorm{\b}^2 $ we get that 
	\begin{align*}
	1 + t^{(k)} (\xk^{\T} \b)  &\geq  1 - t^{(k)} \TNorm{\b}\\
	&\geq 1 - \TNorm{\b}/\TNorm{\b} \geq 0, 
	\end{align*}
	as a result, we get that $ (\vt \b)( \vt \xkp ) \leq 0 $. This holds for an arbitrary eigenvector $\v$ corresponding to the smallest eigenvalue, so $\xkp \in S_E$.
\end{proof}

Thus, in the easy case, it is enough to find some initial vector $\x_0 \in S_E$ to ensure that Riemannian Gradient Descent with the standard metric and projection based retraction, along with step size restriction to $1/\TNorm{\b}$, will converge to a global \btrs{} optimum. Such an $\x_0$ can be trivially found by taking $\x_0 = - \b / \TNorm{\b}$ (if $\b = 0$ we must be in the hard case, and there is no reason to consider $S_E$). 

\begin{rem}
	One might be tempted to forgo the use of the Riemannian gradient in favor of the Euclidean gradient, i.e., use the Projected Gradient Descent iteration 
		\begin{equation}\label{eq:focm.step}
		\xkp = \frac{\xk - t^{(k)} \nabla q (\xk)}{\TNorm{\xk - t^{(k)} \nabla q (\xk}}\ .
		\end{equation}
	Going through the steps of the previous proofs, one can show that this iteration stays in $S_E$ for any step size $t^{(k)}$ when $\matA$ is indefinite, but requires the step size restriction $t^{(k)} < 1/\lmin(\matA)$ when $\matA$ is positive definite. So using iteration~\cref{eq:focm.step} requires spectral information, while the Riemannian iteration~\cref{eq:foro.step} can be used without any knowledge on the spectrum of $\matA$.
\end{rem}

Next, we consider the hard case. According to \Cref{lemma:hardcase_localisglobal_affine_pairs_and_minimizers,lemma:attractors.item.stableglobmin}, in hard cases, the only stable stationary points are global optima. Since we can expect Riemannian Gradient Descent to converge to a stable stationary point for most initial points, we can be tempted to infer that starting the iterations from $\x_0  = -\b/\TNorm{\b}$ will work for the hard case as well. It is however possible that given a carefully crafted starting point Riemannian Gradient Descent will converge to a non optimal stationary point, and indeed for hard case \btrs{} starting from $\x_0  = -\b/\TNorm{\b}$ this is {\em always} the case. Thus, we need a more robust way to select the initial point.  One obvious choice is sampling the starting point from uniform distribution on $\sph$ (or any other continuous distribution on $\sph$). Now we can realistically expect to converge to a stable stationary point. However, we can prove a stronger result. 

First let us define $S_H$ and show that in the hard case, any stationary point in $S_H$ must be a global minimum:
\begin{lemma} \label{lemma:sh.sphere.cond}
	Consider a hard case \btrs{} defined by $\matA$ and $\b$, and let
	$$
	S_H \coloneqq  \{ \x \in \sph~:~ \exists \u \st \matA\u = \lamsub{1}(\matA)\u, \ut \x \neq 0 \}.
	$$
	If $\bar{\x} \in S_H$ is a stationary point, then it is a global optimizer.
\end{lemma}

\begin{proof}
For any eigenvector $ \u $ associated with $ \lamsub{1}(\matA) $, and for any affine eigenpair $ (\bar{\mu}, \bar{\x}) $ we can write:
	$$
	-\ut \b = (\lamsub{1}(\matA) - \bar{\mu})\ut \bar{\x} ~~.
	$$
Without loss of generality, assume that $ \ut \bx \neq 0 $ (since $ \bx \in S_H $ ). 
By our assumption that we are in the hard case, we have that $ \ut \b = 0 $, so it must hold that $ \bar{\mu} = \lamsub{1}(\matA) $. 
Now,~\Cref{lemma:affine_pairs_and_minimizers} guarantees that $\bx$ is a global optimizer.
\end{proof}

Next, we show that if the initial point $\x_0$ is in $S_H$ and the step size is bounded by $1/\TNorm{\b}$, then any accumulation point of the iteration defined by~\Cref{eq:foro.step} must be a global optimum. An equivalent result was established in~\cite[Theorem 4.8]{Beck2018} for \trs{} (that iteration was guaranteed to converge).

\begin{lemma}\label{lemma:unprec.sh.global.convergence}
	Assume that $\b \neq 0$ and that for all eigenvectors $\u$ corresponding to the smallest eigenvalue of $\matA$ we have $\bt \u = 0$. 
	Let $  \{ \xk \}_{k=0}^{\infty}  \subset \sph $ a sequence of iterates obtained by~\Cref{eq:foro.step} with step-size $ t^{(k)} < 1 /  \TNorm{\b} $ for all $ k \in \mathbb{N} $.
	Assume that all accumulation points of $  \{ \xk \}_{k=0}^{\infty}  \subset \sph $ are stationary points. If $\x_0 \in S_H$ then any accumulation point $ \bx $ of $ \{ \xk \}_{k=0}^{\infty} $ is a global minimizer.
\end{lemma}

The proof of~\Cref{lemma:unprec.sh.global.convergence} uses the following auxiliary lemma:
\begin{lemma}\label[lem]{lemma:cost.arq.convergence}
	Let $  \{ \xk \}_{k=0}^{\infty}  \subset \sph $ a sequence of iterates obtained by~\Cref{eq:foro.step}. Assume that for each $k$ we have $q(\xkp) \leq q(\xk)$, and that all accumulation points of $  \{ \xk \}_{k=0}^{\infty}  \subset \sph $ are stationary points. Let $ \bx \in \sph $ an accumulation point of the sequence, then
	\begin{lemlist}

		\item The sequence $  \{ q(\xk) \}_{k=0}^{\infty}$ converges and $\lim_{k \rightarrow \infty} q(\xk) = q(\bx) $. \label{lemma:cost.arq.convergence.cost}
		\item The sequence $  \{ \mu_{\xk} \}_{k=0}^{\infty}$ converges and $ \lim_{k \rightarrow \infty} \mu_{\xk} = \mu_{\bx} $. \label{lemma:cost.arq.convergence.arq}
		\item Let $ \ui $ an eigenvector of $ \matA $ associated with an eigenvalue $ \lamsub{i} $ such that $  \lamsub{i} \neq \mu_{\bx}  $, then the sequence $ \{ |\uit \xk| \}_{k=0}^{\infty} $ converge, and its limit is equal to $  | \uit \b | / | \mu_{\bx} - \lamsub{i} | $. \label{lemma:cost.arq.convergence.uitxk.lim}

	\end{lemlist}
\end{lemma}

\begin{proof}
	$ \bx $ is an accumulation point, so we have a subsequence $  \{ \xki \}_{i=0}^{\infty}  \subset \{ \xk \}_{k=0}^{\infty} $ such that $ \xki \rightarrow \bx $. Since $ \{ q (\xk) \}_{k=0}^{\infty} $ is a monotonic and bounded sequence
	it has a limit. Let us denote the limit $\bar{q}$.
	Thus any subsequence of $ \{ q (\xk) \}_{k=0}^{\infty} $ must converge to $ \bar{q} $, and in particular $ \lim_{i \rightarrow \infty} q(\xki) =\bar{q} $.
	On the other hand, by continuity it holds that $ \lim_{i \rightarrow \infty} q(\xki) = q(\bx) $. Hence $ \bar{q} = q(\bx) $ (i.e., Clause (i) holds.) 
	
	Now suppose in contradiction that $ \mu_{\xk} {\not \rightarrow} \mu_{\bx} $. Then there exists $\delta > 0$ such that for every $ j \in {\mathbb N} $ there is a $ k_j > j $ for which 
	\begin{equation*}
		|\mu_{\xkj} - \mu_{\bx}| >  \delta \label{eq:contra.convergence.arq}.
	\end{equation*}
	Without loss of generality we assume that $\{ \xkj \}_{j=0}^{\infty} $ is convergent, and  denote its limit by $\by$ (since $\{ \xkj \}_{j=0}^{\infty} $ is contained within a compact set, it has a convergent subsequence, and we can chose our sequence to be that subsequence).
	By the lemma's assumptions, both $\bx$ and $ \by $ are stationary points of $q(\cdot)$.
	Since the the map $\x \mapsto \mu_\x$ is continuous we get $ | \mu_{\by} - \mu_{\bx} | \geq \delta $ so clearly $  \mu_{\by} \neq \mu_{\bx} $.
	On the other hand, from the first clause, we have $ q(\by) = q(\bx) $,
	and~\Cref{lemma:uria.1.lucidi} implies that $ \mu_{\by} = \mu_{\bx}  $ arriving at a contradiction, so we must have  $ \mu_{\xk} \rightarrow \mu_{\bx} $ (i.e., Clause (ii) holds.)

	Let ${\cal L}$ be the set of affine eigenvectors corresponding to $\mu_{\bx}$, i.e.,
	$$ 
	{\cal L} \coloneqq \{ \z \in \sph ~|~ (\matA  - \mu_{\bx} \matI_{n} ) \z = - \b \}.  $$ 
	Since $\matA$ is symmetric, every vector $\z \in \sph$ can be decomposed as  $\z = \v + \w$ where $\v$ is in the range of $\matA - \mu_{\bx} \matI_{n}$ and $\w$ is in the null space of $\matA - \mu_{\bx} \matI_{n}$. This implies that for every $\z \in {\cal L}$ we can write $\z = -(\matA - \mu_{\bx} \matI_{n})^+ \b + \w$ where $\w$ is orthogonal to $\z - \w$ and $\z$ has unit norm. 
	Let us define the projection on ${\cal L}$
	$$
	{\cal P}_{\cal L} ( \x ) \coloneqq \argmin_{\z \in {\cal L} } \TNorm{\z - \x}\ ,
	$$
	where in the above ties are broken arbitrarily. 
	 Further define:
	$$
	\rho_\x \coloneqq \x - {\cal P}_{\cal L} ( \x ).
	$$
	Note that $ \TNorm{\rho_{\x}} = \dist (\x , {\cal L}  ) $.
	
	 Considering the above definitions, we write $ | \uit \xk | = | \uit (\rho_{\xk} + {\cal P}_{\cal L} ( \xk )  ) | $ and get the following inequality:
	$$
	| \uit {\cal P}_{\cal L}( \xk ) | - | \uit \rho_{\xk} | \leq | \uit \xk | \leq | \uit {\cal P}_{\cal L}( \xk ) | + | \uit \rho_{\xk} |.
	$$ 
	By construction, $ {\cal P}_{\cal L} ( \xk ) $ is an affine eigenvector corresponding to $\mu_{\bx}  $ for any $ k \in \mathbb{N} $, so recalling that $ \mu_{\bx} \neq \lamsub{i} $ we write:
	\begin{align*}
	\uit {\cal P}_{\cal L} ( \xk ) &= - \uit (( \matA - \mu_{\bx} \matI_{n} )^{ + } \b -  \w), 
	\end{align*}
	where $ \TNorm{( \matA - \mu_{\bx} \matI_{n} )^{ + } \b  - \w} = 1$  (for any matrix $\matX$, $\matX^{+}$ denotes the Moore-Penrose pseudo-inverse of $\matX$) and $ \w $ is in the null space of $ \matA - \mu_{\bx} \matI_{n} $ (if this null-space is empty. Then $ \w = 0 $, and $ \TNorm{( \matA - \mu_{\bx} \matI_{n} )^{ -1 } \b} = 1 $).
	The null space is orthogonal to the range so $  ( \matA - \mu_{\bx} \matI_{n} ) \w = 0 $, which implies that $ \w $ is an eigenvector of $ \matA $ corresponding with the eigenvalue $ \mu_{\bx} \neq \lamsub{i} $ thus $ \uit \w = 0 $ and we have
	\begin{align*}
	\uit {\cal P}_{\cal L} ( \xk ) &= - \uit ( \matA - \mu_{\bx} \matI_{n} )^{ + } \b \\
	&= ( \mu_{\bx} - \lamsub{i} )^{-1} \uit \b .
	\end{align*}
	Now, for any $ \delta > 0 $ there exists a $ K_\delta \in \mathbb{N}$ such that $ \TNorm{\rho_{\xk}} < \delta $ for all $ k \geq K_\delta $, since otherwise we get that there is an accumulation point of $  \{ \xk \}_{k=0}^{\infty}$ outside of $ {\cal L} $, so $\rho_{\xk} \rightarrow 0$. Then, by Cauchy-Schwartz, we have that $ |\uit \rho_{\xk}| \rightarrow 0 $, and  
	\[
	 \lim_{k \rightarrow \infty } | \uit \xk | =  ( \mu_{\bx} - \lamsub{i} )^{-1} | \uit \b | .
	\] 
\end{proof}

\vspace{5mm}

\begin{proof}[Proof of~\cref{lemma:unprec.sh.global.convergence}]
	By~\Cref{eq:foro.step,eq:se.naive.rgrad}, for every $ k=0,1,... $ we have:
	\begin{align}
	\xkp &=\theta_{k}^{1}\xk - \theta_{k}^{2} \grad q(\xk) \nonumber \\
	&=\theta_{k}^{1}\xk - t^{(k)} \theta_{k}^{1} \grad q(\xk) \nonumber \\ 
	&= \theta_{k}^{1}\xk - t^{(k)} \theta_{k}^{1} (\egrad{q}{\xk} - \mu_{\xk} \xk) \nonumber\\
	&= \theta_{k}^{1}(1 + t^{(k)} \mu_{\xk} ) \xk - t^{(k)} \theta_{k}^{1} \egrad{q}{\xk}, \label{eq:theta1k.egradstep.xk}
	\end{align}
	where $ \theta_{k}^{1} \coloneqq \TNorm{\x_{k} - t^{k} \grad q(\x_{k})}^{-1} $ and  $ \theta_{k}^{2} \coloneqq \theta_{k}^{1} t^{(k)} $.
	
	Let $ \bx \in \sph $ an accumulation point of the sequence (thus, by assumption - a stationary point), and $\bar{\mu}$ the associated affine eigenvalue. So:
	\begin{equation}
		-\b = (\matA - \bar{\mu} \matI) \bx . \label{eq:crit.x}
	\end{equation}
	Since $ \x_0 \in S_H$, there exists an eigenvector $ \u $ of $ \matA $ corresponding to $ \lamsub{1} $ such that 
	\begin{equation}
		\ut \x_0 \neq 0. \label{eq:utx0.nonzero}
	\end{equation}
	
	Pre-multiply~\Cref{eq:theta1k.egradstep.xk} by $ \ut $:
	\begin{align}
	\ut \xkp &= \theta_{k}^{1}(1 + t^{(k)} \mu_{\xk} ) \ut \xk - t^{(k)} \theta_{k}^{1} \ut \egrad{q}{\xk} \nonumber \\
	&= \theta_{k}^{1}(1 + t^{(k)} \mu_{\xk} ) \ut \xk - t^{(k)} \theta_{k}^{1} \lamsub{1} \ut \xk \nonumber\\
	&= \theta_{k}^{1}(1 + t^{(k)} (\mu_{\xk} - \lamsub{1}) ) \ut \xk \nonumber  \\ 
	&= \alpha_k  \ut \xk\ , \nonumber 
	\end{align}
	where $ \alpha_k \coloneqq \theta_{k}^{1}(1 + t^{(k)} (\mu_{\xk} - \lamsub{1}) ) $. Note that since:
	\begin{align*}
	\mu_{\xk} - \lamsub{1} &= \xkt \matA \xk + \bt \xk - \lamsub{1} \\
	&\geq \lamsub{1} + \bt \xk - \lamsub{1} \geq - \TNorm{\b}\ ,
	\end{align*}
	and since $ t^{(k)} < \TNorm{\b}^{-1}$, it holds that $ 1 + t^{(k)}(\mu_{\xk} - \lamsub{1}) >0 $ for all $ k $. Combined with the fact that $ \theta_{k}^{1} > 0$, we conclude that $ \alpha_k > 0 $.
	Moreover, for any $ K = 1,2,3... $, it holds that $ |\ut \x_K| = |\ut \x_0 | \prod_{j=0}^{K-1} \alpha_j > 0 $ by~\Cref{eq:utx0.nonzero} (so all iterates are in $S_H$).
	
	Now assume in contradiction that $ \bx $ is not a global minimizer. By~\Cref{lemma:nonglobmin.affine.pairs.minimizers} we have $ \bar{\mu} > \lamsub{1} $.
	Moreover, since $ \b \neq 0 $, there exists an eigenvector $ \v $ of $ \matA  $ associated with an eigenvalue $ \lambda > \lamsub{1} $ such that:
	\begin{equation}
		\vt \b \neq 0. \label{eq:v.def}
	\end{equation}
	Pre-multiplying~\cref{eq:crit.x} by $ \vt $  results in $ (\lambda - \bar{\mu}) \vt \x \neq 0$ hence $ \vt \bx \neq 0 $ and $ \lambda \neq \bar{\mu} $.
	Again, pre-multiply~\Cref{eq:theta1k.egradstep.xk} by $ \vt $:
	\begin{align}
	\vt \xkp &= \theta_{k}^{1}(1 + t^{(k)} \mu_{\xk} ) \vt \xk - t^{(k)} \theta_{k}^{1} \vt \egrad{q}{\xk} \nonumber \\
	&= \theta_{k}^{1}(1 + t^{(k)} (\mu_{\xk} - \lambda) ) \vt \xk - t^{(k)} \theta_{k}^{1} \vt\b . \label{eq:vt.theta1.1} 
	\end{align}
	Since $ \lambda \neq \bar{\mu} $ we have by~\Cref{lemma:cost.arq.convergence.uitxk.lim} that
	\begin{equation}
		 \lim_{k \to \infty} |\vt \xk| = |\vt \b| / |\bar{\mu} - \lambda | > 0, \label{eq:lim.vtxk.nonzero}
	\end{equation}
	so there exists a $ K_1 $ such that $ |\vt \xk| > 0 $ for all $ k \geq K_1 $.
	Hence, for $ k \geq K_1  $ we can write:
	\begin{equation}
	\vt \xkp = \theta_{k}^{1}\left(1 + t^{(k)} \left(\mu_{\xk} - \lambda - \frac{\vt\b}{\vt \xk} \right) \right) \vt \xk\ . \label{eq:vt.theta1.2}
	\end{equation}
	
	Define $ \ressub{\xk} \coloneqq - (\matA - \mu_{\xk} \matI) \xk - \b $, and note that $ \TNorm{\ressub{\xk}} \to 0 $. We have:
	\begin{align*}
		\vt \ressub{\xk} &= - \vt (\matA - \mu_{\xk} \matI) \xk - \vt \b \\
		&= (\mu_{\xk} - \lambda ) \vt \xk - \vt \b,
	\end{align*}
	so, for $ k \geq K_1 $ it holds that 
	\begin{equation}
		\mu_{\xk} - \lambda = \frac{\vt \b + \vt \ressub{\xk}} {\vt \xk}\ . \label{eq:ressub.lambda.muk}
	\end{equation}
	Plugging ~\Cref{eq:ressub.lambda.muk} in \Cref{eq:vt.theta1.2} we get:
	\begin{align*}
	\vt \xkp & = \theta_{k}^{1}\left(1 + t^{(k)} \frac{\vt \ressub{\xk}} {\vt \xk} \right) \vt \xk \\
	         &= \eta_k  \vt \xk\ ,
	\end{align*}
	where $ \eta_k \coloneqq \theta_{k}^{1}(1 + t^{(k)} \vt \ressub{\xk} / \vt \xk ) $. 
	Since $ |\vt \ressub{\xk}| / |\vt \xk| \to 0$ and $ t^{(k)} < \TNorm{\b}^{-1} $, there exists $ K_2 $ such that $ 1 + t^{(k)} \vt \ressub{\xk}/ \vt \xk   > 0$ for all $ k \geq K_2 $, thus $ \eta_k > 0 $ for all $ k \geq K_2 $.
	Consider now the difference $ \alpha_k - \eta_k $:
	\[
	(\alpha_k - \eta_k) / \theta_{k}^{1} = t^{(k)} \left(\mu_{\xk} - \lamsub{1} - \frac{\vt \ressub{\xk}} {\vt \xk} \right),
	\]
	and again, since $ |\vt \ressub{\xk}| / |\vt \xk| \to 0 $ by construction of $ \ressub{\xk} $ and~\Cref{lemma:cost.arq.convergence.uitxk.lim}, combined with that $ \mu_{\xk} - \lamsub{1} \to \bar{\mu} - \lamsub{1} > 0 $ by~\Cref{lemma:cost.arq.convergence.uitxk.lim,lemma:hardcase_localisglobal_affine_pairs_and_minimizers}, there is a $ K_3 $ such that $ \alpha_k > \eta_k $ for all $ k \geq K_3 $.
	
	Let $ K \coloneqq \max\{K_1 , K_2, K_3 \} $, then for any $ k \geq K $ it holds that:
	\[
	0 < \eta_k \leq \alpha_k\ .
	\]
	Moreover, let $ k \geq K $ then:
	\begin{align*}
		|\vt \xkp| &= |\vt \x_K| \prod_{j=K}^{k} \eta_j \\
		&\leq |\vt \x_K| \prod_{j=K}^{k} \alpha_j \\
		&=  \frac{|\vt \x_K|}{|\ut \x_K|} |\ut \x_K| \prod_{j=K}^{k} \alpha_j \\
		&=  \frac{|\vt \x_K|}{|\ut \x_K|}  |\ut \xkp|\ ,
	\end{align*}
	where we used the fact  that $ |\ut \x_K| \neq 0 $ for all $K$.
	Taking the limit $ k \to \infty $, we have that:
	\begin{align*}
		\lim_{k \to \infty} |\vt \xk| &\leq \frac{|\vt \x_K|}{|\ut \x_K|} \lim_{k \to \infty} |\ut \xk|  = 0 .
	\end{align*}
	in contradiction to~\cref{eq:lim.vtxk.nonzero}. Thus $ \bar{\x} $ must be a global minimizer.
\end{proof}

We are left with the task of choosing an initial point $ \x_0 \in S_H$. Unlike the case of choosing a point in $S_E$, we cannot devise a deterministic method for choosing $\x_0 \in S_H$ without spectral information (the eigenvectors corresponding to the minimal eigenvalues of $\matA$). Nevertheless, it is possible to choose a random initial point which almost surely is in $S_H$ as long as $\matA$ is not a scalar-matrix, i.e., a multiple of $\matI_{n}$ (if $\matA$ is a scalar-matrix, solving \btrs{} is trivial). When $\matA$ is not a scalar-matrix, the set $\sph \smallsetminus S_H$ has measure 0. So choosing an initial point by sampling any continuous distribution on $\sph$ (e.g., Haar measure) will almost surely be in $S_H$.

Without spectral information, it is impossible to know a priori if we are dealing with the easy or the hard case. 
So, we employ the double start idea suggested by Beck and Vaisbourd~\cite{Beck2018}: we execute two iterations, one in $S_E$ and the other in $S_H$, and on conclusion the solution with minimum objective value is chosen. This summarized in \Cref{alg:naive.convergent}.

\begin{algorithm}[H]
	\caption{Double start, Riemannian optimization algorithm for globally solving \btrs{}.}\label{alg:naive.convergent}
	\begin{algorithmic}[1]
		\Input  $\matA \in \RR^{n\times n}$, $\b \in \RR^{n}$ 
		\State$ \x_0 \gets -\b / \TNorm{\b} $
		\label{startingpoint_se}
		\State $\bar{\x},\bar{\mu}\gets$\Call{NaiveRGD}{$\matA$,$\b$,$\x_0$}~(\cref{alg:rgd.canonical}) \label{line:alg.naive.fs} 
		\State {Sample $ \x_0 $ from some continuous distribution on $ \sph $}
		\State $\tilde{\x},\tilde{\mu}\gets$\Call{NaiveRGD}{$\matA$,$\b$,$\x_0$}~(\cref{alg:rgd.canonical}) \label{line:alg.naive.ds}
		\State $ \xopt \gets \argmin_{} \{ q (\z) ~|~ \z \in \{ \tilde{\x},  \bar{\x} \} \} $
		\State \Return $ \xopt, \muopt $
	\end{algorithmic}
\end{algorithm}

\begin{algorithm}[H]
	\caption{Riemannian Gradient Descent subroutine for solving \btrs{} using standard geometry.}\label{alg:rgd.canonical}
	\begin{algorithmic}[1]
		\Function{NaiveRGD}{$\matA$, $\b$, $\x_0$} 
		\Params  $ \tau \in (0,1),c \in (0,1) $ 
		\For{$ k = 0,1,2,... $} 
		\State $\eta_k \gets \grad q(\x_k) = \proj{\x_k}\egrad{q}{\x_k} = (\matI_{n} - \x_k \x_k^\T)(\matA \x_k + \b)$ 
		\State $ t^{(k)} \gets 1 / \TNorm{\b} $
		\While{$ q(\x_k) - q(R_{\x_k} (-t^{(k)} \eta_k)) < t^{(k)} c\TNormS{\eta_k }  $}
		\State $ t^{(k)} \gets \tau t^{(k)} $
		\EndWhile
		\State $ \x_{k+1} \gets R_{\x_k}(- t^{(k)} \eta_k )$
		\EndFor
		\State \Return $ \x_k, \mu_{(\x_k)} $
		\EndFunction
	\end{algorithmic}
\end{algorithm}

	\section{Generic Riemannian \btrs{} Solver}\label{sec:solving.prec}

In the previous section we presented a method that globally solves \btrs{} and uses Riemannian Gradient Descent with the standard Riemannian metric and projection based retraction. It is desirable to lift these restrictions, and in particular the requirements to use only Riemannian Gradient Descent and the standard Riemannian metric. It is well known that typically Riemannian Conjugate Gradient  enjoys faster convergence rates compared to Riemannian Gradient Descent, and that incorporating a non-standard Riemannian metric, a technique termed {\em Riemannian preconditioning} in the literature, may introduce considerable acceleration~\cite{MS16,Shustin2019}. In this section we propose such an algorithm. However, unlike \Cref{alg:naive.convergent}, the algorithm presented in this section requires spectral information: an eigenvector corresponding to the minimal eigenvalue of $ \matA $. 

As already mentioned, while we can expect Riemannian optimization algorithms to converge to a stable stationary point, such points might be local minimizers that are not global. We need to detect whether convergence to such a point has occurred, and somehow handle this. The key observation is the following lemma, which shows that if we have sufficiently converged (see the condition on the residual in \Cref{lemma:analysis.bound.progress}) to any stationary point other than the global solution we can devise a new iterate which reduces the objective. Furthermore, the reduction in the objective function is bounded from below, and this guarantees that it will be executed a finite number of times. The lemma is useful only to easy case \btrs{} (in hard cases the global minimums are the only stable stationary points).

\begin{lemma}
	\label{lemma:analysis.bound.progress}
	Suppose that $\u$ is an eigenvector corresponding to the minimal eigenvalue of $\matA$ such that $\ut \b \neq 0$. Let $\alpha = | \b^{\T} \u |$. Consider some candidate approximate affine eigenvector $\x\in\sph$ and its corresponding affine Rayleigh quotient $\mu_\x$, and suppose that $\mu_\x > \lamsub{1}(\matA)$. Let 
	$\rb_\x = \mu_\x \x  - \matA \x - \b$ (i.e., $\rb_\x$ is the residual in upholding the affine eigenvalue equation). Let 
	$$\xlpr \coloneqq \x - 2 (\ut \x) \u.$$
	If  $\TNorm{\rb_\x} \leq \alpha / 2$ then 
	$$
	q(\x) - q(\xlpr) \geq \frac{\alpha^2}{(\mu_{\max} - \lamsub{1}(\matA))}\ ,
	$$
	where $\mu_{\max}$ is the maximal affine eigenvalue of $(\matA, \b)$.	
\end{lemma}

\begin{rem}
The transformation $\x \mapsto \xlpr$ is a special case of a more general transformation suggested in~\cite{Lucidi1998}, and the initials LPR in the subscript correspond to the authors name. 
\end{rem}

\begin{proof}
	\label{lemma:analysis.bound.progress.proof}
	Note that $ q(\xlpr) = q(\x) - 2(\u^{\T} \b )(\u^{\T} \x) $.
	Pre-multiplication of $ \mu_\x \x  $ by $ \ut $ enables to write 
	\begin{eqnarray*}
	\mu_\x \u^{\T} \x &=& \u^{\T} \matA \x + \u^{\T} \b + \u^{\T} \rb_\x \\
	&=& \lamsub{1} \u^{\T}\x + \u^{\T} \b + \u^{\T} \rb_\x\ .
	\end{eqnarray*}
	And thus 
	\begin{equation}
	\u^{\T} \x = \frac{ \u^{\T} \b + \u^{\T} \rb_\x}{\mu_\x -  \lamsub{1} }\ ,
	\end{equation}
	where we used the fact that $\mu_\x \neq  \lamsub{1} $.
	Multiplying by $ 2 (\u^{\T} \b ) $ results in 
	\begin{eqnarray*}
	2 (\u^{\T} \b )(\u^{\T} \x )  &=& 2  \frac{ (\u^{\T} \b )^2 + (\u^{\T} \rb_\x)(\u^{\T} \b )}{\mu_\x - \lamsub{1}}\\
	&\geq& 2 \frac{ \alpha^2 - \TNorm{\rb_\x} \alpha }{\mu_\x - \lamsub{1}} \\
	&\geq& 2 \frac{ \alpha^2 -  \alpha^2 / 2 }{\mu - \lamsub{1}} \\
	&\geq& \alpha^2 / (\mu_{\max} - \lamsub{1}),
	\end{eqnarray*}
	thus $ q(\xlpr) = q(\x) - 2(\u^{\T} \b )(\u^{\T} \x) \leq q(\x) - \alpha^2 / (\mu_{\max} - \lamsub{1}) $.
\end{proof}

We can leverage this observation in the following way. First, we use an eigensolver to find an eigenvector $\u$ corresponding to the minimal eigenvalue of $\matA$ such that $\ut\b \neq 0$. If no such eigenvector exists, then we are in the hard case, and we use a Riemannian optimization solver with initial point $\x_0$ sampled from the Haar measure on $\sph$. The only stable stationary points are the global optimizers, so we expect the solver to converge to a global optimizer. If, however, we found such a vector $\u$, we are in the easy case, and there might be a stable  stationary point other than the global minimizer. 

Now, we use an underlying Riemannian optimization solver, augmenting its convergence test with the requirement that $\TNorm{\rb} \leq \alpha/2$ where $\alpha = |\ut \b|$. Once the Riemannian optimization solver returns $\x$, we check whether $\mu_\x < \lamsub{1}(\matA)$. If it is, then we return $\x$. Otherwise, we replace $\x$ with $\xlpr$ and restart the Riemannian optimization. The algorithm is summarized in \Cref{alg:prec}. We have the following theorem:

\begin{theorem}
 	Consider executing \Cref{alg:prec} on an easy case \btrs{}. Then:
 	\begin{thmlist}
 		\item If the algorithm terminates, it returns a point $\bx$ such that $\mu_{\bx} < \lamsub{1}(\matA)$.
 		\item If all intermediate Riemannian iterations (line \ref{line:run-R}) reduce the objective value, then the number of occasions where line \ref{line:lpr} is visited is finite, and the algorithm will terminate in finite time. 
 	\end{thmlist}
\end{theorem}
 \begin{spacing}{1.01}
\begin{algorithm}[H]
    
 \caption{Generic Riemannian solver for \btrs{}.}   
 
 \begin{algorithmic}[1]
    \label{alg:prec}
  \INPUT  $\matA\in\R^{n\times n}$, $\b\in\R^n$, and underlying Riemannian solver ${\cal R}$.

  \State Sample $ \xp $ from $ \mathcal{N}(0,1)^n $ and set $ \x \gets \xp / \TNorm{\xp}$. \label{alg:prec.starting_point.start}
  \State \label{line:eig} Use an eigensolver to find an eigenvector $\u$ corresponding to the minimal eigenvalue of $\matA$ for which $\ut\b \neq 0$.
  \State If none such exist: return the result of ${\cal R}$ starting from $\x$.

  \If{$ (\u^{\T} \b)(\u^{\T} \x) > 0 $}
  	\State $ \x\gets -\x $
  \EndIf \label{alg:prec.starting_point.end}
  \State $\alpha \gets |\ut \b|$ 
  \Loop \label{alg:prec.while.start}
  	\State \label{line:run-R} Run ${\cal R}$, augmenting its convergence criteria with the additional criteria that \State $\TNorm{\rb_\x}\leq \alpha/2$, starting from $\x$ to obtain $\bx$.
  	\If{$ \mu_{\bar{\x}} < \lamsub{1}(\matA) $}
	  	\State \label{line:ret}\Return $\bar{\x}$
	\Else
		\State \label{line:lpr} $\x \gets \bx_{\noun{LPR}}$
	\EndIf \label{alg:prec.lucidi.end}
  \EndLoop \label{alg:prec.while.end}
  	
 \end{algorithmic}
\end{algorithm}
\end{spacing}

\begin{proof}
	Since we are in the easy case, a $\u$ will be found in line~\ref{line:eig}. 
	So the algorithm may return only via line~\ref{line:ret}. 
	This requires that $\mu_{\bx} < \lamsub{1}(\matA)$ as required.
	Since intermediate applications of ${\cal R}$ reduce the objective, and the transformation in line \ref{line:lpr} reduces the objective, the objective is always decreased. 
	Furthermore, since line \ref{line:lpr} is executed only when $\TNorm{\rb_\x} < \alpha/2$, \Cref{lemma:analysis.bound.progress} guarantees that the reduction in the objective that occurs in line \ref{line:lpr} is lower bounded by a constant, so the amount of such reductions is finite, and so is the number of times line~\ref{line:lpr} is executed.   
\end{proof}
 
Although the algorithm allows for running the underlying solver ${\cal R}$ multiple times, we expect that in non-pathological cases it will execute at most twice. The reason is that there are at most two stable stationary points. If we sufficiently converge to a stable stationary point other than the global solution, we expect line~\ref{line:lpr} to push the objective below the objective value of the local minimizer, and future descents will be towards to global minimum.

	\section{Preconditioned Solver}\label{subsec:precond-scheme}
		\Cref{alg:prec} uses an underlying Riemannian solver ${\cal R}$. The running time of \Cref{alg:prec} is highly dependent on how fast ${\cal R}$ converges to a stationary point. In this section we propose a framework for incorporating a preconditioner into ${\cal R}$ for the case that ${\cal R}$ is some standard general-purpose Riemannian algorithm (e.g., Riemannian Gradient Descent and Riemannian Conjugate Gradients). The idea is to use Riemannian preconditioning. That is, the preconditioner is incorporated by using a non-standard Riemannian metric on $\sph$.

To construct our framework, We first define a smooth mapping $\x\mapsto\matM_\x$ from the sphere to the set of symmetric positive definite matrices. We now endow $\RR^n$ with the metric $\bar{g}_{\x}(\eta_\x,\xi_\x) = \eta^{\T}_\x \matM_\x \xi_\x$. We then consider $\sph$ as a Riemannian submanifold of $\RR^n$ endowed with this metric, thus the Riemannian metric on $\sph$ is defined by
$g_{\x}(\eta_\x,\xi_\x) = \eta^{\T}_\x \matM_\x \xi_\x$ in ambient coordinates. We refer to the mapping $\x\mapsto\matM_\x$ as the {\em preconditioning scheme}.

Recall that the analysis in \Cref{sec:solving.prec} was independent of the choice of metric.
In this section we analyze how the preconditioning scheme affects convergence rate, and use these insights to propose a preconditioning scheme based on a constant seed preconditioner $\matM$.

To study the effect of the preconditioning scheme on the rate of convergence, we analyze the spectrum of the Riemannian Hessian at stationary points, when viewed as a linear operator on the tangent space.
Such analyses are well motivated by the literature, see \cite[Theorem 4.5.6, Theorem 7.4.11 and Equation 7.50]{Absil2008}, though these results are, unfortunately, only asymptotic.
The following theorem provides bounds on the extreme eigenvalues of the Riemannian Hessian at stationary points.

\begin{theorem}\label{thm:cond-hess}
	Suppose that $\bx$ is a stationary point of \btrs{}, and that $\bar{\mu}$ is the corresponding affine eigenvalue.
	Then, the spectrum of $\Hess q(\bx)$ is contained in the interval
	\[
	\left[
		\lambda_{\min} \left(
			\matM_{\bx}^{-1/2}\left[\matA - \bar{\mu}\matI_{n}\right]\matM_{\bx}^{-1/2}
		\right),
		\lambda_{\max} \left(
			\matM_{\bx}^{-1/2} \left[\matA-\bar{\mu}\matI_{n}\right] \matM_{\bx}^{-1/2}
		\right)
	\right].
	\]
	Furthermore, if $\x_\star$ is a global optimum, and we are in the easy case, then
	\[
		\kappa(\Hess q(\x_\star))\leq\kappa\left(\matM_{\x_{\star}}^{-1/2}\left[\matA - \mu_{\star}\matI_{n}\right]\matM_{\x_{\star}}^{-1/2}\right),
	\]
	where $\kappa(\cdot)$ denotes the condition number of a matrix.
\end{theorem}

\begin{proof}\label{thm:cond-hess.proof}
	\Cref{tab:formulas} gives a formula, in ambient coordinates,
	for the Riemannian Hessian at a stationary $\bx$; for all $\eta_{\bx}\in \tsp{\bx}$ we have
	\[
	\Hess q(\bx)[\eta_{\bx}]=\proj{\bx}\matM_{\bx}^{-1}\left[\matA-\bar{\mu}\matI_{n}\right]\eta_{\bx}\,.
	\]
	As the Riemannian Hessian operator is self-adjoint with respect to the Riemannian
	metric~\cite[Proposition 5.5.3]{Absil2008}, it is possible to use the Courant-Fisher Theorem to get that for every point $\z\in\sph$
	\[
	\lambda_{\max}(\Hess q(\z))=\max_{0\neq\eta_{\z}\in \tsp{\z} }\frac{g_{\z}(\eta_{\z},\Hess q(\z)[\eta_{\z}])}{g_{\z}(\eta_{\z},\eta_{\z})}\ ,
	\]
	\[
	\lambda_{\min}(\Hess q(\z))=\min_{0\neq\eta_{\z}\in \tsp{\z} }\frac{g_{\z}(\eta_{\z},\Hess q(\z)[\eta_{\z}])}{g_{\z}(\eta_{\z},\eta_{\z})}\ .
	\]
	The above is stated in a coordinate-free manner. In ambient coordinates, viewing $\tsp{\bx}$ as a subspace of $\RR^{n}$, and using the specific formula for the Hessian at stationary points, we have:
	\begin{align*}
		\lambda_{\max}(\Hess q(\bx)) &= \max_{\eta_{\bx}\neq 0, \, \eta_{\bx}^{\T} \bx =0} \frac{\eta_{\bx}^{\T} \matM_{\bx}\proj{\bx}\matM_{\bx}^{-1}\left[\matA-\bar{\mu}\matI_{n}\right]\eta_{\bx}}{\eta_{\bx}^{\T}\matM_{\bx}\eta_{\bx}} \\
		& = \max_{\eta_{\bx}\neq 0, \, \eta_{\bx}^{\T} \bx =0} \frac{\eta_{\bx}^{\T}\matP_{\bx}^{\T}\left[\matA-\bar{\mu}\matI_{n}\right]\eta_{\bx}} {\eta_{\bx}^{\T}\matM_{\bx}\eta_{\bx}} \\
		& = \max_{\eta_{\bx}\neq 0, \, \eta_{\bx}^{\T} \bx =0} \frac{\eta_{\bx}^{\T}\left[\matA-\bar{\mu}\matI_{n}\right]\eta_{\bx}} {\eta_{\bx}^{\T}\matM_{\bx}\eta_{\bx}} \\
		& = \max_{\xi_{\bx}\neq0,\,\xi_{\bx}^{\T}\matM_{\bx}^{-1/2}\bx=0} \frac{\xi_{\bx}^{\T}\matM_{\bx}^{-1/2}\left[\matA-\bar{\mu}\matI_{n}\right]\matM_{\bx}^{-1/2}\xi_{\bx}}{\xi_{\bx}^{\T}\xi_{\bx}}\\
		& \leq \lambda_{\max}\left(\matM_{\bx}^{-1/2}\left[\matA-\bar{\mu}\matI_{n}\right]\matM_{\bx}^{-1/2}\right).
	\end{align*}
	where the second equality we used the fact that $\matM_{\y}\matP_{\y}\matM_{\y}^{-1}=\matP_{\y}^{\T}$ (see~\cite[Section E.2]{Shustin2019}). For the third equality we used the fact that $\matP_{\bx}$ is a projector on $\tsp{\bx}$, thus for any $\eta_{\bx} \in \tsp{\bx} $,  we have $ \proj{\bx} \eta_{\bx}=\eta_{\bx}$.
	For the last inequality we used the Courant-Fisher Theorem yet again.

	Similarly,
	\[
	\lambda_{\min}(\Hess q(\bx)) \geq \lambda_{\min} \left(\matM_{\bx}^{-1/2}\left[\matA-\bar{\mu}\matI_{n} \right]\matM_{\bx}^{-1/2}\right).
	\]
	This proves the first part of the theorem.

	As for the second part, since we are in the easy case and $\x_\star$ is the global \btrs{} minimizer, then $\mu_{\x_{\star}} < \lambda_{\min}$.
	So, $\matA-\mu_{\x_{\star}}\matI_{n}$ is positive definite, and
	\[
	\kappa(\Hess q(\bx)) \leq \kappa\left(\matM_{\x_\star}^{-1/2}\left[\matA-\mu_{\x_\star}\matI_{n}\right]\matM_{\x_\star}^{-1/2}\right).
	\]
\end{proof}

\Cref{thm:cond-hess} presents a new perspective on the "hardness" of the hard case: suppose that we are in the hard case, and that $\xopt$ is a global optimum with a corresponding affine eigenvalue $ \muopt $.  Then $\muopt=\lamsub{1}$ and $\matA - \muopt \matI_{n}$ is singular - a case for which our bounds are meaningless (note that the theorem requires the easy case settings). If, however, we approach the hard case in the limit, then the condition number explodes. This result echoes the analysis of Carmon and Duchi, in that \btrs{} instances can get arbitrarily close to being a hard-case, requiring more iterations in order to find an exact solution~\cite{Carmon2018}.

We now leverage \Cref{thm:cond-hess} to propose a systematic way to build a preconditioning scheme using some seed preconditioner $\matM$. In light of \Cref{thm:cond-hess}, one would want to have $\matM_{\xopt}\approx\matA - \muopt \matI_{n}$ where $\xopt$ denotes the global optimum.   We approximate each term separately.
First, we replace $\matA$ with some approximation $\matM$ of $\matA$ forming the \textit{seed preconditioner}.
Next we approximate the second term $-\muopt \matI_{n}$.
Exact calculation of this quantity raises a fundamental problem as it requires us to know $\xopt$, which after all, is the vector we are looking for, so  $\muopt$ remains inaccessible and this is where the varying metric comes into play.
Given a good approximation $\x \in \sph$ of $\xopt$, we know that $\mu_{\x}$ is a good guess for $\muopt$, i.e, $\mu_{\x} = \argmin_{\mu}\TNormS{ \left( \matA-\mu \matI_{n} \right)\x+\b}$.
Furthermore, it is easy to see that $\mu_{\xk} \to \muopt $ when $\xk \to \xopt $.
Therefore, we can approximate $-\muopt \matI_{n}$ with $-\mu_{\x} \matI_{n}$.
The matrix $\matM_\x$ is then formed by adding these two approximations, but with an additional filter applied to $-\mu_{\x}$, resulting in:
\begin{equation}\label{eq:precond-scheme}
\matM_{\x}=\matM + \phi(-\mu_{\x})\matI_{n}\ ,
\end{equation}
where $\phi(\cdot)$ is a smooth function such that for all $\alpha\in\RR$ we have:
\phantomsection{\label{item:phi.requirement}}
\begin{enumerate}
	\item $\phi(\alpha) > - \lmin(\matM)$ - making sure that $ \matM + \phi(\alpha) \matI_{n} $ is  positive definite.\label{item:phi.requirement.posdef}
	\item $ \phi(\alpha) \approx \max(-\lmin(\matM),\alpha) $ - so that  $ \matA - \muopt \matI_{n}$  is well approximated by $ \matM + \phi(-\muopt) \matI_{n} $  for values of $ \x $ near the global solution $ \xopt $. \label{item:phi.requirement.goodapprox}
\end{enumerate}
The filter $\phi(\cdot)$ is designed so that $\matM_\x$ is always positive definite and that the mapping $\x \mapsto \matM_\x$ is smooth.
One possible concrete construction of $\phi$ is detailed in~\Cref{subsec:phi}.
Note that the preconditioning scheme shown in~\Cref{eq:precond-scheme} is fully defined by $\matM$.
Hence, we refer to $\matM$ as the \emph{seed preconditioner} (or just \emph{preconditioner}).

Obviously, the incorporation of our preconditioning scheme defined by~\Cref{eq:precond-scheme} in Riemannian optimization algorithms requires solving one or more linear equations, whose matrix is $\matM_{\x}$ with each iteration in order to compute the Riemannian gradient.
Recalling that $\matM_{\x}$ is a scalar-matrix shift of the seed-preconditioner $\matM$, we also need $\matM$ that is amenable to the solution of such systems (e.g., $\matM$ is low rank and we use the Woodbury formula).


%
%

	\section{Solving the \trs{}\label{sec:trs}}
\newcommand{\hatAdef}{\ensuremath{\left[\begin{matrix}
			0 & \\
			& \matA
		\end{matrix}\right]}}
\newcommand{\hatbdef}{\ensuremath{\left[
		\begin{matrix}
			0 \\ \b
		\end{matrix}
		\right]}}
\newcommand{\zopt}{\ensuremath{\z_\star}}
\newcommand{\slice}[3]{\ensuremath{{#1}_{\text{#2:#3}}}}

In previous sections we presented algorithms for finding a global solutino of a \btrs{} by means of Riemannian optimization. We now consider the \trs{}, and show how we can leverage our algorithm for solving \btrs{} to solve \trs{}.

It is quite common to encounter $\matA$ and $\b$ such that the global solution of \trs{} and \btrs{} coincide. In fact, they will always coincide if $\matA$ is not positive definite. When $\matA$ is positive definite, the global solution of the $\trs{}$ is equal to $\matA^{-1} \b$ if that vector is in the interior of $\BB^n$ (and so in this case \btrs{} and \trs{} solutions do not coincide). If that vector is not in the interior, there exists a solution on the boundary, and the solution of \trs{} and \btrs{} do coincide. In light of these observations, a trivial algorithm for solving \trs{} via \btrs{} is given in \Cref{alg:naive.trs.via.btrs}.

\begin{algorithm}[H]
	\begin{algorithmic}[1]
		\INPUT $ \matA, \b $
		\If{$ \lamsub{1}(\matA) > 0 $}
			\State $ \Delta \gets \TNormS{\matA^{-1}\b} $
			\If{$ \Delta < 1 $}
				\State \Return $\matA^{-1} \b $ 
			\Else
				\State solve \btrs{}
			\EndIf
		\Else
			\State solve \btrs{}
		\EndIf
	\end{algorithmic}
	\caption{\label{alg:naive.trs.via.btrs}Solving \trs{} via \btrs{} }
\end{algorithm}

One obvious disadvantage of \Cref{alg:naive.trs.via.btrs} is that it requires us to either have a priori knowledge of whether $\matA$ is positive definite, or to somehow glean whether $\matA$ is positive definite or not (e.g., by testing positive definiteness~\cite{HanEtAl17,Bakshi2020TestingPS}).
This method might also require to compute $ \matA^{-1}\b$, which can be costly as well. 
An alternative method, and one that is superior when $\matA$ is positive definite, is to use the augmentation trick\footnote{In spite of the method's  simplicity, we are not aware of any descriptions of this method earlier then Phan et al.'s (relatively) recent work.} suggested by Phan et al.~\cite{Phan2020}. 
The augmentation trick is as follows. 
Given a \trs{} defined by $\matA$ and $\b$, construct an augmented \btrs{}  by
\begin{equation} \label{eq:aug.A.hat}
\begin{aligned}
    \hat{\matA} &=\hatAdef \in \RR^{(n+1) \times (n+1)},  \\
    \hat{\b} &= \hatbdef \in \RR^{n+1}. 
\end{aligned}
\end{equation}

One can easily see that a solution of \trs{} can be obtained by discarding the first coordinate of the augmented \btrs{}  solution. 

If $\matA$ is positive definite, the augmented \btrs{} is necessarily in the hard case. One can see that directly, but is also evident from the fact that there are at least two global solutions, one obtained from the other by flipping the sign of first coordinate (an easy case \btrs{} has only one global solution). If $\matA$ is not positive definite, the augmented \btrs{} is a hard case if and only if \trs{} is a hard case \trs{}.

In case $\matA  $ is positive semidefinite, the augmented \btrs{} may be of either the hard or easy case.
The next lemma shows that for this specific case, we can deterministically find an  initial vector in the intersection of $S_H$ and $S_E$, thus making the double-start strategy redundant.

\begin{lemma}\label{lemma:spd_single_init}
	Let $\hatA $ and $ \bhat $ be constructed from $\matA$ and $\b$ according to~\cref{eq:aug.A.hat}. Assume that $ \matA $ is symmetric positive semidefinite and define 
	\begin{equation}
		\x_0 = \left[ \begin{matrix} 1/\sqrt{\TNormS{\b} + 1} \\ - \b / \sqrt{\TNormS{\b} + 1} \end{matrix} \right]\ , \label{eq:global_x0_spd_cannonical}
	\end{equation}
	then $\x_0 \in S_H \cap S_E$.
\end{lemma}

\begin{proof}
	First, note that indeed $ \x_0 \in \sphn$.
	If $\matA$ is strictly positive definite, then the minimal eigenvalue of $\hatA$ is $0$, and the only unit-norm eigenvectors corresponding to it are $\pm\e_1$ where $\e_1 $ is a unit vector with 1 on its first coordinate and the remaining entries are zeros. Since $\bhat^\T \e_1 = 0$, we have that $S_E = \sphn$ and so $\x_0 \in S_E$. Since $\x^\T_0 \e_1 = 1/\sqrt{\TNorm{\b}^2 + 1} \neq 0$, we get that $\x_0 \in S_H$.
	
	Next, consider the case that $\matA$ is positive semidefinite but not positive definite, i.e., it is singular. Now, in addition to $\pm \e_1$, there are additional eigenvectors that correspond to the $0$ eigenvalue. The conditions for inclusion in $S_E$ and $S_H$ based on $\pm \e_1$ were already verified, so we focus on the rest of the eigenvalues. We need to consider only eigenvectors that are orthogonal to $\e_1$. Such an eigenvector must have the structure $\hat{\u}=[0; \u]$ where $\matA \u = 0$. Then, 
		\begin{equation*}
		(\bthat \hat{\u}) (\x_0^{\T} \hat{\u}) = (\bt \u)\left(- \bt \u / \sqrt{\TNormS{\b} + 1} \right) \leq 0.
		\end{equation*}
		So $\x_0 \in S_E$.  In addition, inclusion in $S_H$ still holds since $\e_1$ is still an eigenvector corresponding to the minimal eigenvalue (that is 0), and we already argued that $\x^\T_0 \e_1 \neq 0$.
\end{proof}

%


	\section{Numerical Illustrations}\label{sec:experiments}
		We illustrate \Cref{alg:naive.convergent,alg:prec} on three synthetically generated sets of matrices. One corresponds to an easy case \btrs{}, the second to a hard case \btrs{}, and the third, while technically an easy case, is ``almost hard". 

The method for generating test matrices is based on the method in~\cite{Beck2018}, adding the slight modification of defining the spectrum of $ \matA $ as a mixture of equispaced ``signal'' and random ``noise''. Namely, a random symmetric matrix $ \matA $ of dimension $ n = 2000 $ is generated, where 75\%   of $ \matA $'s eigenvalues are sampled from a normal distribution with zero mean and standard deviation of $10^{-3}$. The rest of $ \matA $'s spectrum is equispaced in \text{[-5,10]}.
The expected level of difficulty of each problem is determined by the gap between  $ \mu_{\star} $ and the eigenvalue of $\matA$ that is closest to $\mu_{\star}$, that is $\lmin(\matA)$.  
The gap is set to $2,~ 10^{-8},~0$ to simulate the easy, almost hard, and the hard case respectively.  $\xopt $ is sampled at random from $ \sph $.
Once $ \matA , ~ \mu_{\star}$ and $ \xopt $ are set,  $ \b  $ is obtained by solving $ (\matA - \mu_{\star} \matI_{n} ) \mat{y} =  - \xopt $ for $\mat{y}$. 
For each difficulty level (determined by $\lmin (\matA) - \mu_{\star} $) we produce 20 synthetic instances of the \btrs{} as described above.

We use \Cref{alg:prec} with a preconditioner. To build the seed preconditioner, we use a fixed-rank symmetric sketch similar to the method presented in ~\cite{Tropp2016}, to get a symmetric matrix $ \matM $ of rank 50 and  its spectral decomposition.

\begin{figure}[H]
	\centering
	\includegraphics[width=\textwidth]{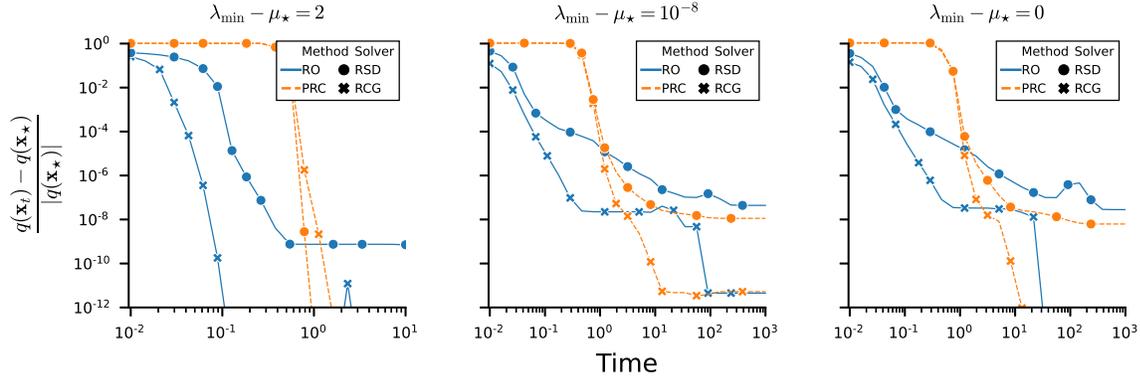}
	
	\caption{The progression in \textbf{objective} (relative) error where x axis shows the \textbf{time} in seconds.\label{fig:t_obj}}
\end{figure}

\begin{figure}[H]
	\centering
	\includegraphics[width=\textwidth]{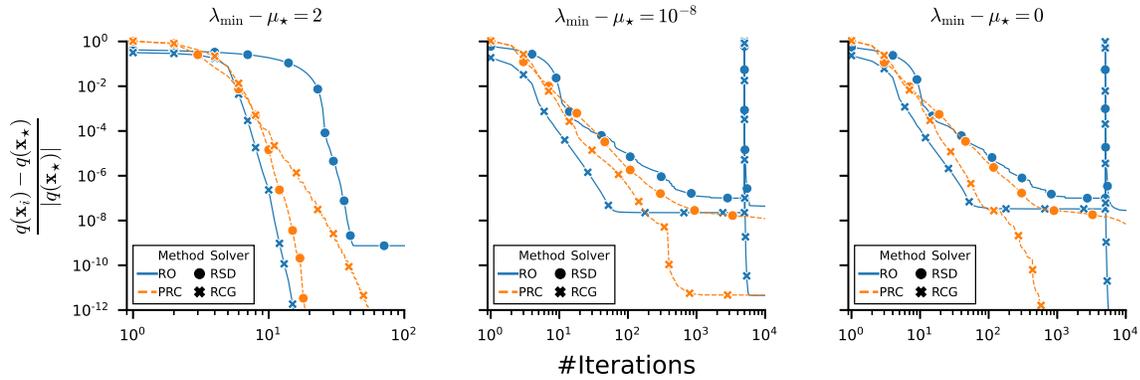}
	
	\caption{The progression in \textbf{objective} (relative) error where x axis shows the \textbf{iteration number}.\label{fig:it_obj}}
\end{figure}

\begin{figure}[H]
	\centering
	\includegraphics[width=\textwidth]{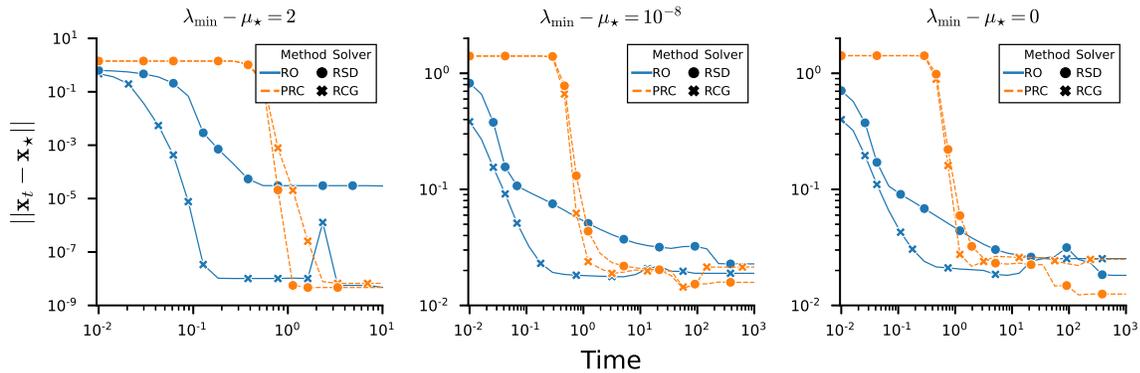}
	
	\caption{The progression in \textbf{argument} error  where x axis shows the \textbf{time} in seconds.\label{fig:t_arg}}
\end{figure}

\begin{figure}[H]
	\centering
	\includegraphics[width=\textwidth]{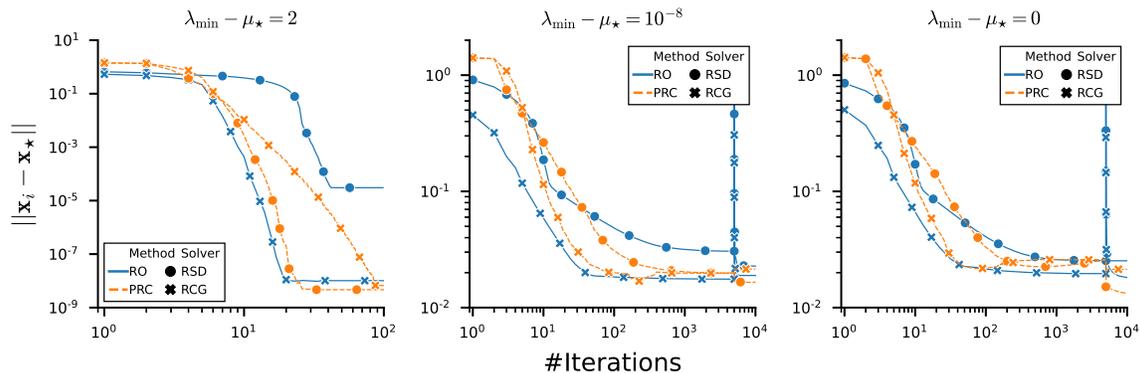}
	
	\caption{The progression in \textbf{argument} error where x axis shows the \textbf{iteration number}.\label{fig:it_arg}}
\end{figure}

Results are reported in~\Cref{fig:it_arg,fig:it_obj,fig:t_arg,fig:t_obj}. Plain Riemannian optimization is labeled as ``RO'', while preconditioned Riemannian optimization is labeled as "PRC". For each approach, we used both Riemannian Steepest Descent (RSD) and Riemannian Conjugated Gradients (RCG) solvers. In general, our limited experiments suggest that RCG does a much better job than RSD. When the problem is very well-conditioned (i.e., it is an easy case), RCG does a much better job than it's preconditioned counterpart. This is due to the preprocessing cost of the preconditioned approach, and is not uncommon when using preconditioning for well-conditioned problems. In contrast, for the hard case and almost-hard case, we see a clear benefit for preconditioned CG. With respect to RSD, preconditioning almost always help. 

When considering the progression in terms of number of iterations, $\TNorm{\x_{i} - \xopt}$ and $|(q(\x_{i}) - q(\xopt))/q(\xopt)|$ for the hard and almost hard cases, it is apparent that in the 5000'th iterations, the ``RO'' run instances suffer a `bump' in their values. 
This bump caused by line~\ref{line:alg.naive.ds} of ~\Cref{alg:naive.convergent} that forces a re-start of the optimization process from a new, random starting point after the first iteration process in line~\ref{line:alg.naive.fs} which is set to finish when convergence criteria $\TNorm{\grad q(\x_{i})} \leq 10^{-12}$ is reached or following 5000 iterations.
The phenomena is not observed for the preconditioned~\Cref{alg:prec}, since re-initialization within the loop specified by line~\ref{line:run-R} are made from points with values lower than those of previous iteration (hence the objective is ever decreasing).

One final remark is in order. Note that when the problem is hard or almost-hard, even though the algorithms finds a near-minimizer, the argument error, $\|\x_{t} - \x_{\star}\|$, is large. This is expected given our observation that hardness of \btrs{} translates to ill conditioning (of the Riemannian Hessian).

\newpage

	\bibliographystyle{plain}
    \bibliography{bibfile}

	\newpage
	\appendix
	\section{Constructing \texorpdfstring{$ \phi $}{phi}}\label{subsec:phi}
\DeclareRobustCommand{\cnameref}[1]{%
  \namecref{#1} \nameref{#1}%
}%
\DeclareRobustCommand{\Cnameref}[1]{%
  \nameCref{#1} \nameref{#1}%
}

We now show a simple way to construct a smooth $\phi(\cdot)$ fulfilling both requirements stated in~\cref{item:phi.requirement.posdef,item:phi.requirement.goodapprox} in~\Cref{item:phi.requirement}. 

First, we choose some time parameter $\epsilon > 0$ and set $ d \coloneqq \lmin(\matM) - \epsilon $. We construct $\phi(\alpha)$ to be a smoothed out over-estimation of $ \max(\alpha, -d) $. We first construct an under estimation. Let $\gamma>1$ be another parameter, and define
\[
\varphi(\alpha)=\frac{\alpha+d}{2}\left(1-\tanh\left(-\gamma\left(\alpha+d\right)\right)\right)-d.
\]
Then $ \varphi(\alpha) $ is a smooth function that approximates $ \max(\alpha, -d)  $, but it is an under estimation: $\varphi(\alpha) \leq  \max(\alpha, -d) $. 

\begin{figure}[H]
		\centering
		\input{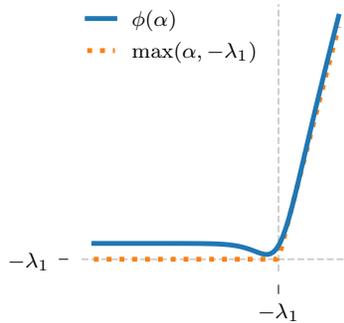}
		\caption{Illustration of the approximation's behavior. Values of the variable $ \alpha $ are depicted in the x axis. Note that  $ \phi(-\alpha) > - \lmin(\matM) $ for all $ \alpha $, while in addition $ \phi(\alpha) $ well approximates $ \alpha $ for values of $ \alpha \geq -  \lmin(\matM) $.}
		\label{fig:phi}
\end{figure}

While the difference between  $ \max(\alpha, -d)  $ and $\varphi(\alpha)$ reduces significantly when $ \gamma $ grows, we want to make sure that our approximation is greater than or equal to the $ \max(\alpha, -d) $. To that end, let 
$$
\alpha_0 \coloneqq - \frac{\operatorname{W}\left[0,e^{-1}\right]+1}{2\gamma} - d, 
$$
where $ \operatorname{W}\left[0,\cdot \right] $ is the zero branch of the Lambert-$ \operatorname{W}$ function. Now, set:
\begin{equation}
\phi(\alpha)\coloneqq\varphi(\alpha)-\varphi(\alpha_{0})-d.\label{eq:phi}
\end{equation}
See~\Cref{fig:phi} for a graphical illustration of $\phi$.
It is possible to show that 
\begin{equation}
\label{eq:errphi}
0 \leq \phi(\alpha) -  \max(\alpha, -\lmin(\matM)) \leq \frac{\operatorname{W}\left[0,e^{-1}\right]+1}{2\gamma}(1 - \tanh((\operatorname{W}\left[0,e^{-1}\right]+1)/2)) )   + \epsilon,
\end{equation}
so~\Cref{item:phi.requirement.posdef} holds, and the approximation error (\Cref{item:phi.requirement.goodapprox}) is small if $\epsilon$ is sufficiently small and $\gamma$ is sufficiently large. The proof of \Cref{eq:errphi} is rather technical and does not convey any additional insight on the \btrs{} and its solution, so we omit it.

\end{document}